\documentclass[aos,MSNbibl,seceqn,nameyear,rotating,dvips]{arximspdf}

\doi{10.1214/14-AOS1219} 
\volume{42}
\issue{4}
\pubyear{2014}
\firstpage{1233}
\lastpage{1261}

\makeatletter
\newcommand{\rrvert}{\vert}
\newcommand{\llvert}{\vert}

\newtheorem{lem}{Lemma}
\newcommand{\argmax}{\mathop{\arg\max}}
\newtheorem{theorem}{Theorem}
\newproclaim{remark}{Remark}
\makeatother

\begin{document}
\begin{frontmatter}

\title{A second-order efficient empirical Bayes confidence~interval}
\runtitle{An empirical Bayes confidence interval}

\begin{aug}
\author[a]{\fnms{Masayo}~\snm{Yoshimori}\corref{}\thanksref{t2}\ead[label=e1]{masayo@sigmath.es.osaka-u.ac.jp}}
\and
\author[b]{\fnms{Partha}~\snm{Lahiri}\thanksref{t3}\ead[label=e2]{plahiri@survey.umd.edu}}
\runauthor{M. Yoshimori and P. Lahiri}
\thankstext{t2}{Supported by the JSPS KAKENHI Grant Number 242742.}
\thankstext{t3}{Supported by the NSF SES-085100.}
\affiliation{Osaka University and University of Maryland}
\address[a]{Department of Medical Innovation\\
Osaka University Hospital\\
2-2 Yamadaoka\\
Suita, Osaka 5650871\\
Japan\\
\printead{e1}}

\address[b]{The Joint Program in Survey Methodology\\
University of Maryland\\
1218 Lefrak Hall\\
College Park, Maryland 20742\\
USA\\
\printead{e2}}
\end{aug}

\received{\smonth{3} \syear{2013}}
\revised{\smonth{3} \syear{2014}}

%
\begin{abstract}
We introduce a new adjusted residual maximum likelihood method (REML)
in the context of producing an empirical Bayes (EB) confidence interval
for a normal mean, a problem of great interest in different small area
applications. Like other rival empirical Bayes confidence intervals
such as the well-known parametric bootstrap empirical Bayes method, the
proposed interval is second-order correct, that is, the proposed
interval has a coverage error of order $O(m^{-{3}/{2}})$. Moreover,
the proposed interval is carefully constructed so that it always
produces an interval shorter than the corresponding direct confidence
interval, a property not analytically proved for other competing
methods that have the same coverage error of order $O(m^{-{3}/{2}})$. The proposed method is not simulation-based and requires
only a fraction of computing time needed for the corresponding
parametric bootstrap empirical Bayes confidence interval. A Monte Carlo
simulation study demonstrates the superiority of the proposed method
over other competing methods.
\end{abstract}

%
\end{frontmatter}

\section{Introduction}\label{sec1}

Fay and Herriot (\citeyear{FH}) considered
empirical Bayes estimation of small
area means $\theta_i$ using the following two-level Bayesian model and
demonstrated, using real life data, that they outperform both the
direct and synthetic (e.g., regression) estimators.

\textit{The Fay--Herriot model}:

 For $i=1,\ldots,m$,
\begin{longlist}
\item[Level 1] (sampling distribution): $y_i|\theta_i
\stackrel{\mathrm{ind}}{\sim}
N(\theta_i,D_i)$;
\item[Level 2] (prior distribution): $\theta_i
\stackrel{\mathrm{ind}}{\sim}
N(x_i^{\prime}\beta, A)$.
\end{longlist}
In the above model, level 1 is used to account for the sampling
distribution of the direct survey estimates $y_i$, which are usually
weighted averages of the sample observations in area $i$. Level 2 prior
distribution links the true small area means $\theta_i$ to a vector of
$p<m$ known area level auxiliary variables $x_i=(x_{i1},\ldots,
x_{ip})^{\prime}$, often obtained from various administrative records.
The hyperparameters $\beta\in R^p$, the $p$-dimensional Euclidean
space, and $A\in[0,\infty)$ of the linking model are generally unknown
and are estimated from the available data.

It is often difficult or even impossible to retrieve all important
sample data within small areas due to confidentiality or other reasons
and the only data an analyst may have access to are aggregate data at
the small area level. The Fay--Herriot model comes handy in such
situations since only area level aggregate data are needed to implement
the model. Even when unit level data are available within small areas,
analysts may have some preference for the Fay--Herriot model over a
more detailed (and perhaps more scientific) unit level model in order
to simplify the modeling task. One good feature of the Fay--Herriot
model is that the resulting empirical Bayes (EB) estimators of small
area means are design-consistent. In the Fay--Herriot model, sampling
variances $D_i$ are assumed to be known, which often follows from the
asymptotic variances of transformed direct estimates [Efron and Morris
(\citeyear{EM75}), \citet{CR}] and/or from empirical variance modeling
[\citet{FH}]. This known sampling variance assumption causes
underestimation of the mean squared error (MSE) of the resulting
empirical Bayes estimator of $\theta_i$. Despite this limitation, the
Fay--Herriot model has been widely used in different small area
applications [see, e.g., \citet{CR}, Efron and Morris (\citeyear{EM75}),
\citet{FH,Bell.2007}, and others].

Note that the empirical Bayes estimator of $\theta_i$ obtained by \citet{FH} can be motivated as an empirical best prediction
(EBP) estimator [in this case same as the empirical best linear
unbiased prediction (EBLUP) estimator] of the mixed effect $\theta
_i=x_i^{\prime}\beta+v_i$, under the following linear mixed model:
\[
y_i=\theta_i+e_i=x_i^{\prime}
\beta+v_i+e_i, \qquad i=1,\ldots,m,
\]
where the $v_i$'s and $e_i$'s are independent with $v_i \stackrel
{\mathrm{i.i.d.}}{\sim} N(0,A)$ and $e_i \stackrel{\mathrm{ind}}{\sim}
N(0, D_i)$; see \citet{PR.1990} and \citet{Rao.2003}.

In this paper, we consider interval estimation of small area means
$\theta_i$. An interval, denoted by $I_i$, is called a $100(1-\alpha)\%$
interval for $\theta_i$ if $P(\theta_i \in I_i |\beta,A)
=1-\alpha$, for any fixed $\beta\in R^p, A\in(0,\infty)$, where the
probability $P$ is with
respect to the Fay--Herriot model. Throughout the paper, $P(\theta_i
\in I_i |\beta,A)$ is referred to as the coverage
probability of the interval $I_i$; that is, coverage is defined in
terms of the joint distribution of $y$ and $\theta$ with fixed
hyperparameters $\beta$ and $A$.
Most intervals proposed in the literature can be written as: $\hat\theta
_i\pm s_\alpha\hat\tau_i(\hat\theta_i)$, where $\hat\theta_i$ is an
estimator of $\theta_i$, $\hat\tau_i(\hat\theta_i)$ is an estimate of
the measure of uncertainty of $\hat\theta_i$ and $s_\alpha$ is suitably
chosen in an effort
to attain coverage probability close to the nominal level $1-\alpha$.

Researchers have considered different choices for $\hat\theta_i$. For
example, the choice $\hat\theta_i=y_i$ leads to the direct confidence
interval $I_i^{D}$, given by
\[
I_i^{D}\dvtx y_i\pm z_{\alpha/2}
\sqrt{D_i},
\]
where $z_{\alpha/2}$ is the upper $100(1-\alpha/2)\%$ point of $N(0,1)$.
Obviously, for this direct interval, the coverage probability is $1
-\alpha$. However, when $D_i$ is large as in the case of small area
estimation, its length is too large
to make any reasonable conclusion.

The choice $\hat\theta_i=x_i^{\prime}\hat\beta$, where $\hat\beta$ is a
consistent estimator of $\beta$, provides an interval based on the
regression synthetic estimator of $\theta_i$. \citet{HM}
considered this choice with $\hat\tau_i(\hat\theta_i)=\sqrt{\hat A}$,
$\hat A$ being a consistent estimator of $A$, and obtained $s_\alpha$
using a parametric bootstrap
method. This approach could be useful when $y_i$ is missing for the
$i$th area.

We call an interval empirical Bayes (EB) confidence interval if we
choose an empirical Bayes estimator for $\hat\theta_i$. There has been
a considerable interest in constructing empirical Bayes confidence
intervals, starting from the work of \citet{Cox} and \citet{Ma},
because of good theoretical and empirical properties of empirical Bayes
point estimators. Before introducing an empirical Bayes confidence
interval, we introduce the Bayesian credible interval in the context of
the Fay--Herriot model.
When the hyperparameters $\beta$ and $A$ are known, the Bayesian
credible interval
of $\theta_i$ is obtained using the posterior distribution of $\theta
_i\dvtx \theta_i|y_i;(\beta,A)\sim N[\hat{\theta}_i^{B},
\sigma_i(A)]$, where $\hat\theta_i^{\mathrm B}\equiv\hat{\theta
}_i^{B}(\beta, A) =(1-B_i)y_i+B_ix_i^{\prime}\beta, B_i\equiv
B_i(A)=\frac{D_i}{D_i+A}, \sigma_i(A)=\sqrt{\frac{AD_i}{A+D_i}}\ (i=1,\ldots,m)$. Such a credible interval is
given by
\[
I_i^{B}(\beta,A)\dvtx \hat{\theta}_i^{B}(
\beta,A)\pm z_{\alpha/2} \sigma_i (A).
\]
The Bayesian credible interval cuts down the length of the direct
confidence interval by $100\times(1-\sqrt{1-B_i})\%$ while
maintaining the exact coverage $1-\alpha$ with respect to the joint
distribution of $y_i$ and $\theta_i$. The maximum benefit from the
Bayesian methodology is achieved
when $B_i$ is close to 1, that is, when the prior variance $A$ is much
smaller than the sampling variances $D_i$.

In practice, the hyperparameters are unknown. \citet{Cox} initiated the
idea of developing an one-sided empirical Bayes confidence interval for
$\theta_i$ for a special case of the Fay--Herriot model with $p=1,
x_i^{\prime}\beta=\beta$ and $D_i=D\ (i=1,\ldots,m)$.
The two-sided version of his confidence interval is given by
\[
I_i^{\mathrm{Cox}}(\hat\beta,\hat A_{\mathrm{ANOVA}})\dvtx \hat{
\theta}_i^{B}(\hat\beta,\hat A_{\mathrm{ANOVA}})\pm
z_{\alpha/2} \sigma(\hat A_{\mathrm{ANOVA}}),
\]
where $\hat{\theta}_i^{B}(\hat\beta,\hat A_{\mathrm{ANOVA}})=(1-\hat B)y_i+\hat
B\hat\beta$, an empirical Bayes estimator of $\theta_i$;
$\hat\beta=m^{-1}\sum_{i=1}^my_i$ and $\hat B=D/(D+\hat A_{\mathrm{ANOVA}})$ with
$\hat A_{\mathrm{ANOVA}}=\operatorname{max} \{(m-1)^{-1}\sum_{i=1}^m(y_i-\hat\beta
)^2-D, 0 \}$. An extension of this ANOVA
estimator for the Fay--Herriot model can be found in \citet{PR.1990}.

Like the Bayesian credible interval, the length of the Cox interval is
smaller than that
of the direct interval. However, the Cox empirical Bayes confidence
interval introduces
a coverage error of the
order $O(m^{-1})$, not accurate enough in most small area applications.
In fact, \citet{Cox} recognized
the problem and considered a different $\alpha^{\prime}$, motivated
from a higher-order asymptotic expansion, in order
to bring the coverage error down to $o(m^{-1})$. However, such an
adjustment may cause the interval to be undefined when
$\hat A_{\mathrm{ANOVA}}=0$ and sacrifices an appealing feature of
$I_i^{\mathrm{Cox}}(\hat\mu,\hat A_{\mathrm{ANOVA}})$, that is, the length of such
interval may no longer be less than that of the direct method.

One may argue that Cox's method has an undercoverage problem because it
does not incorporate uncertainty due to estimation of the regression
coefficients $\beta$
and prior variance $A$ in measuring uncertainty of the empirical Bayes
estimator of $\theta_i$. \citet{Ma} used an improved measure of
uncertainty for his empirical Bayes estimator that incorporates the
additional uncertainty due to the estimation of the model parameters.
Similar ideas can be found in \citet{PR.1990} for a more general model.
However, Basu, Ghosh and Mukerjee (\citeyear{Basu}) showed that the coverage error of the
empirical Bayes confidence interval proposed by \citet{Ma} remains
$O(m^{-1})$. In the context of the Fay--Herriot model, Diao et al.
(\citeyear{Diao}) examined the higher order asymptotic coverage of a class of
empirical Bayes confidence intervals of the form: $\theta_i^{\mathrm{EB}}\pm
z_{\alpha/2}\sqrt{\operatorname{mse}_i}$, where $\theta_i^{\mathrm{EB}}$ is an empirical
Bayes estimator of $\theta_i$ that uses a consistent estimator of $A$
and $\operatorname{mse}_i$ is a second-order unbiased estimator of $\operatorname
{MSE}(\theta_i^{\mathrm{EB}})$ given in Datta and Lahiri (\citeyear{DL00}). They showed
that the coverage error for such an interval is $O(m^{-1})$. In a
simulation study, Yoshimori (\citeyear{Y}) observed poor finite sample
performance of such empirical Bayes confidence intervals. Furthermore,
it is not clear if the length of such confidence interval is always
less than that of the direct method.
\citet{Mb} considered a variation of his (\citeyear{Ma}) empirical Bayes
confidence interval where he used a hierarchical Bayes-type point
estimator in place of the previously used empirical Bayes estimator and
conjectured, with some evidence, that the
coverage probability for his interval is at least $1-\alpha$. He also
noted that the coverage probability tends to
$1-\alpha$ as $m$ goes to $\infty$ or $D$ goes to zero. However,
higher-order asymptotic properties of this confidence interval are unknown.

Using a Taylor series expansion, Basu, Ghosh and Mukerjee (\citeyear{Basu}) obtained
expressions for the order $O(m^{-1})$ term of the coverage errors of
the Morris' interval and another prediction interval proposed by Carlin
and Louis [(\citeyear{CL}), page~98], which were then used to calibrate the lengths
of these empirical Bayes confidence intervals in order to reduce the
coverage errors down to $o(m^{-1})$.
However, it is not known if the lengths of their
confidence intervals are always smaller than that of the direct method.
Using a multilevel model,
\citet{Nan} obtained an empirical Bayes confidence interval for a
small area mean and showed that
asymptotically it converges to the nominal coverage probability.
However, he did not study the higher-order
asymptotic properties of his interval.

Researchers considered improving the coverage property of the Cox-type
empirical Bayes confidence interval by changing
the normal percentile point $z_{\alpha/2}$. For the model used by \citet{Cox}, \citet{LL} proposed a
prediction interval based on parametric bootstrap samples. However, the
order of their coverage error has not been studied analytically.
\citet{DGSL} used a Taylor series approach similar to that of
Basu, Ghosh and Mukerjee (\citeyear{Basu})
in order to calibrate the Cox-type empirical Bayes confidence
interval for the general Fay--Herriot model. Using mathematical tools
similar to Sasase and \citet{SK}, Yoshimori (\citeyear{Y}) extended the
method of \citet{DGSL} and Basu, Ghosh and Mukerjee (\citeyear{Basu}) when REML
estimator of $A$ is used.

For a general linear mixed model, Chatterjee, Lahiri and Li (\citeyear{CLL.2008}) developed a
parametric bootstrap empirical Bayes confidence interval for a general
mixed effect and examined its higher order asymptotic properties. For
the special case, this can be viewed as a Cox-type empirical Bayes
confidence interval where $z_{\alpha/2}$ is replaced by percentile
points obtained using a parametric bootstrap method. While the
parametric bootstrap empirical Bayes confidence interval of Chatterjee,
Lahiri and Li (\citeyear{CLL.2008}) has good theoretical properties,
one must apply caution
in choosing~$B$, the number of bootstrap replications, and the
estimator of~$A$. In two different simulation studies, \citet{L.L.2010} and Yoshimori (\citeyear{Y}) found that the parametric bootstrap
empirical Bayes confidence interval did not perform well when REML
method is used to estimate $A$. \citet{L.L.2010} developed an
adjusted REML estimator of $A$ that works better than the REML in their
simulation setting. Moreover, in absence of a sophisticated software,
analysts with modest computing skills may find it a daunting task to
evaluate parametric bootstrap confidence intervals in a large scale
simulation experiment.
The coverage errors of confidence intervals developed by \citet{DGSL}, Chatterjee, Lahiri and Li (\citeyear{CLL.2008}) and \citet{L.L.2010} are of the
order $O(m^{-{3}/{2}})$. However, there is no analytical result
that suggests the lengths of these
confidence intervals are smaller than the length of the direct method.

In Section~\ref{sec2}, we introduce a list of notation and regularity
conditions used in the paper.
In this paper, our goal is to find an empirical Bayes confidence
interval of $\theta_i$ that (i) matches the coverage error properties
of the best known empirical Bayes method such as the one proposed by
Chatterjee, Lahiri and Li (\citeyear{CLL.2008}), (ii) has length smaller than that of the
direct method and (iii) does not rely on simulation-based heavy
computation. In Section~\ref{sec3}, we propose such a new interval method for
the general Fay--Herriot model by replacing the ANOVA estimator of $A$
in the Cox interval by a carefully devised adjusted residual maximum
likelihood estimator of $A$. \citet{L.L.2009} introduced a
generalized (or adjusted) maximum likelihood method for estimating
variance components in a general linear mixed model. \citet{L.L.2010} and \citet{YL} examined different adjustment
factors for point estimation of the small area means in the context of
the Fay--Herriot model. But none of the authors explored adjusted
residual likelihood method for constructing small area confidence
intervals. In Section~\ref{sec4}, we compare our proposed confidence interval
methods with the direct, different Cox-type EB confidence intervals and
the parametric bootstrap empirical Bayes confidence interval method of
Chatterjee, Lahiri and Li (\citeyear{CLL.2008}) using a Monte Carlo simulation study. The
proofs of all technical results presented in Section~\ref{sec3} are deferred to
the \hyperref[app]{Appendix}.

\section{A list of notation and regularity conditions}\label{sec2}
We use the following notation throughout the paper:
\begin{description}
\item$y=(y_1,\ldots,y_m)^{\prime}$, a $m\times1$ column vector of
direct estimates;
\item$X^{\prime}=(x_1,\ldots,x_m)$, a $p\times m$ known matrix of rank $p$;
\item$q_i=x_{i}^{\prime}(X^{\prime}X)^{-1}x_{i}$, leverage of area $i$
for level 2 model, $(i=1,\ldots,m)$;
\item$V=\operatorname{diag}(A+D_1,\ldots,A+D_m)$, a $m\times m$ diagonal matrix;
\item $P=V^{-1}-V^{-1}X(X^{\prime}V^{-1}X)^{-1}X^{\prime}V^{-1}$;
\item$L_{\mathrm{RE}}(A)=|X'V^{-1}X|^{-{1}/{2}}|V|^{-{1}/{2}}\exp
(-\frac{1}{2}y^{\prime}Py )$,
the residual likelihood function of $A$;
\item$h_i(A)$ is a general area specific adjustment factor;
\item$L_{i;\mathrm{ad}}(A)= h_i(A)\times L_{\mathrm{RE}}(A)$, adjusted residual
likelihood function of $A$ with a general adjustment factor $h_i(A)$;
\item$\hat A_{h_i}=\argmax_{A\in[0,\infty)}L_{i;\mathrm{ad}}(A)$, adjusted
residual maximum likelihood estimator of $A$ with respect to a general
adjustment factor $h_i(A)$;
\item$l_{\mathrm{RE}}(A)=\log[L_{\mathrm{RE}}(A)]$;
\item$\tilde{l}_{i;\mathrm{ad}}(A)=\log h_i(A)$;
\item$l_{i;\mathrm{ad}}(A)=\log L_{i;\mathrm{ad}}(A)$;
\item$\tilde{l}_{i,\mathrm{ad}}^{(k)}(A) \equiv\frac{\partial^k \tilde
{l}_{i,\mathrm{ad}}(A)}{\partial A^k}$, $k$th derivative of $\tilde
{l}_{i,\mathrm{ad}}(A), (k\ge1)$;

\item${l}_{i,\mathrm{ad}}^{(k)}(A) \equiv\frac{\partial^k
{l}_{i,\mathrm{ad}}(A)}{\partial A^k}$, $k$th derivative of ${l}_{i,\mathrm{ad}}(A),
(k\ge1)$;

\item$\hat V=\operatorname{diag}(\hat A_{h_1} +D_1,\ldots, \hat A_{h_m} +D_m),
(i=1,\ldots,m)$;
\item$\tilde{\beta}=(X^{\prime}{V}^{-1}X)^{-1}X^{\prime}{V}^{-1}y$,
weighted least square estimator of $\beta$ when $A$ is known;
\item$\hat\beta=(X^{\prime}{\hat V}^{-1}X)^{-1}X^{\prime}{\hat
V}^{-1}y$, weighted least square estimator of $\beta$ when $A+D_i$ is
replaced by $\hat A_{h_i}+D_i, (i=1,\ldots,m)$;
\item${B}_i=D_i/(\hat A +D_i)$, shrinkage factor for the $i$ area,
$(i=1,\ldots,m)$;
\item$\hat{B}_i\equiv\hat B_i(\hat A_{h_i})=D_i/(\hat A_{h_i}
+D_i)$, estimated shrinkage factor for the $i$ area, $(i=1,\ldots,m)$;
\item $\hat\theta_i^{\mathrm B}\equiv\hat{\theta}_i^{B}(\beta, A)
=(1-B_i)y_i+B_ix_i^{\prime}\beta$;
\item$\hat{\theta} _i^\mathrm{EB} \equiv\hat{\theta} _i^\mathrm{EB}
(\hat A_{h_i})\equiv\hat{\theta} _i^\mathrm{B} (\hat\beta, \hat
A_{h_i})=(1-\hat{B}_i)y_i + \hat{B}_i x_i'\hat{\beta}$, empirical Bayes
estimator of $\theta_i, (i=1,\ldots,m)$;

\item$I_i^{\mathrm{Cox}}( \hat\beta, \hat A_{h_i})\equiv I_i^{\mathrm{Cox}}(\hat
A_{h_i})\dvtx \hat{\theta}_i^{B}(\hat\beta, \hat A_{h_i})\pm z_{\alpha/2}
\sigma_i (\hat A_{h_i})$, Cox-type EB confidence interval of $\theta_i$
using adjusted REML $\hat A_{h_i}$, where $z=z_{\alpha/2}$ is the upper
$100(1-\alpha/2)\%$ point of the normal deviate.
\end{description}
We use the following regularity conditions in proving different results
presented in this paper.

\textit{Regularity conditions}:
\begin{description}
\item{R1}: The logarithm of the adjustment term $\tilde
{l}_{\mathrm{ad}}(A)$ [or $\tilde{l}_{i,\mathrm{ad}}(A)$] is free of $y$ and is five
times continuously differentiable with respect to $A$. Moreover, the
$g$th power of the $|\tilde{l}_{\mathrm{ad}}^{(j)}(A)|$ [or $|\tilde
{l}_{i,\mathrm{ad}}^{(j)}(A)|$] is bounded for $g > 0$ and $j=1,2,3,4,5$;
\item{R2}: $\operatorname{rank}(X)=p$;
\item{R3}: The elements of $X$ are uniformly bounded implying
$\sup_{j\geq1}q_j=O(m^{-1})$;
\item{R4}: $0<\inf_{j\geq1}D_j\leq\sup_{j\geq1}D_j<\infty$,
$A\in(0,\infty)$;
\item{R5}: $|\hat{A}_{h_i}|<C_{+}m^{\lambda}$, where $C_{+}$ a
generic positive constant and $\lambda$ is small positive constant.
\end{description}

\section{A new second-order efficient empirical Bayes confidence interval}\label{sec3}
We call an empirical Bayes interval of $\theta_i$ second-order
efficient if the coverage error is of order $O(m^{-{3}/{2}})$ and
length shorter than that of the direct confidence interval. The goal of
this section is to produce such an interval that requires a fraction of
computer time required by the recently proposed parametric bootstrap
empirical Bayes confidence interval. Our idea is simple and involves
replacement of the ANOVA estimator of $A$ in the empirical Bayes
interval proposed by \citet{Cox} by a carefully devised adjusted
residual maximum likelihood estimator of $A$.

Theorem~\ref{th1} provides a higher-order asymptotic expansion of the
confidence interval $I_i^{\mathrm{Cox}}( \hat A_{h_i})$. The theorem holds for
any area $1\le i\le m$, for large $m$.

\begin{theorem}\label{th1}
Under regularity conditions \textup{R1--R5}, we have
%
\begin{equation}
P \bigl\{\theta_i \in I_i^{\mathrm{Cox}}( \hat
A_{h_i}) \bigr\}=1-\alpha+z\phi (z)\frac{a_i+b_i[h_i(A)]}{m}+O
\bigl(m^{-{3}/{2}}\bigr),\label{PI.cox}
\end{equation}
where
%
\begin{eqnarray}\label{ai}
a_i&=&-\frac{m}{\operatorname{tr}(V^{-2})} \biggl[\frac{4D_i}{A(A+D_i)^2}+\frac
{(1+z^2)D_i^2}{2A^2(A+D_i)^2}
\biggr]
\nonumber
\\[-8pt]
\\[-8pt]
\nonumber
&&{}- \frac{mD_i}{A(A+D_i)}x_{i}^{\prime}\operatorname{Var}(\tilde{
\beta})x_i,
\\
b_i&\equiv& b_i\bigl[h_i(A)\bigr]=
\frac{2m}{\operatorname{tr}(V^{-2})}\frac{D_i}{A(A+D_i)}\times \tilde{l}^{(1)}_{i;\mathrm{ad}}.\label{bi}
\end{eqnarray}
\end{theorem}

We can produce higher order asymptotic expansion of the coverage
probability of Cox-type EB confidence interval with any standard
likelihood-based estimator of $A$ available in the literature (e.g.,
residual maximum likelihood, profile maximum likelihood, different
adjusted residual and profile maximum likelihood, etc.) simply by
choosing an appropriate $h_i(A)$ [e.g., for REML, $h_i(A)=1$] and using
equation (\ref{PI.cox}). We have verified that coverage errors for all
these Cox-type EB confidence intervals are of order $O(m^{-1})$. We
can, however, use equation (\ref{PI.cox}) to reduce the coverage error
to the order $O(m^{-3/2})$ by choosing $h_i(A)$ such that the order
$O(m^{-1})$ term in the right-hand side of (\ref{PI.cox}) vanishes.
More specifically, we first obtain an expression for $h_i(A)$ by
finding a solution to the following differential equation:
%
\begin{equation}
a_i+b_i \bigl[h_i(A)
\bigr]=0\label{diffeqn}
\end{equation}
and then maximize the adjusted residual likelihood $L_{i;\mathrm{ad}}(A)$ with
respect to $A\in[0,\infty)$ to obtain our adjusted residual maximum
likelihood estimator of $A$, which is used to construct the desired
Cox-type second-order efficient EB confidence interval for $\theta_i$.
Notice that we can produce two different new adjusted REML estimators
of $A$ by using generalized least square (GLS) and ordinary least
square (OLS) estimators of $\beta$ in the EB estimator of $\theta_i$.
Let $h_{i;{\mathrm{gls}}}(A)$ and $h_{i;{\mathrm{ols}}}(A)$ denote the adjustment factors
that are solutions of $h_i(A)$ in (\ref{diffeqn}) with GLS and OLS
estimators of $\beta$ in $\hat\theta_i^{\mathrm{EB}}$, respectively. We denote
the corresponding adjusted residual maximum likelihood estimators of
$A$ by $\hat{A}_{i;{\mathrm{gls}}}$ and $\hat{A}_{i;{\mathrm{ols}}}$. Note that in general we
cannot obtain $h_{i;{\mathrm{ols}}}(A)$ as a special case of $h_{i;{\mathrm{gls}}}(A)$ except
for the balanced case $D_i=D, i=1,\ldots,m$ when the GLS and OLS
estimators of $\beta$ are identical. Consequently, $\hat A_{i;{\mathrm{gls}}}$ is
generally different from $\hat A_{i;{\mathrm{ols}}}$ except for the balanced case
when $\hat A_{i;{\mathrm{gls}}}=\hat A_{i;{\mathrm{ols}}}=\hat A_i$
(say).

Theorem~\ref{th2} provides expressions for $h_{i;{\mathrm{gls}}}(A)$ and $h_{i;{\mathrm{ols}}}(A)$
and states the uniqueness of $\hat A_i$ for the balanced case. In
Theorem~\ref{th2} and elsewhere in the paper, $C$ is a generic constant free of $A$.

\begin{theorem}\label{th2}
\textup{(i)} The expressions for $h_{i;{\mathrm{gls}}}(A)$ and $h_{i;{\mathrm{ols}}}(A)$ are
given by
%
\begin{eqnarray}
h_{i;{\mathrm{gls}}}(A)&=&CA^{({1}/{4})(1+z^2)}(A+D_i)^{({1}/{4})(7-z^2)}\label{h.unbalance.gls}
\nonumber
\\[-8pt]
\\[-8pt]
\nonumber
&&{}\times\exp\biggl[\int\frac{1}{2}\operatorname{tr}\bigl(V^{-2}\bigr)
x_i^{\prime}\bigl(X^{\prime
}V^{-1}X
\bigr)^{-1}x_i\,dA\biggr],
\\
h_{i;{\mathrm{ols}}}(A)&=&CA^{({1}/{4})(1+z^2)}(A+D_i)^{({1}/{4})(7-z^2)}
\Biggl[\prod_{i=1}^{m}(A+D_i)
\Biggr]^{({1}/{2})q_i}\label
{h.unbalance.ols}
\nonumber
\\[-8pt]
\\[-8pt]
\nonumber
&&{}\times\exp \biggl[-\frac{1}{2}\operatorname{tr}\bigl(V^{-1}\bigr)
x_i^{\prime}\bigl(X^{\prime
}X\bigr)^{-1}X^{\prime}VX
\bigl(X^{\prime}X\bigr)^{-1}x_i \biggr].
\end{eqnarray}
\textup{(ii)} For the balanced case $D_i=D$ ($i=1, \ldots, m$), we have
%
\begin{equation}
h_{i;{\mathrm{gls}}}(A)=h_{i;{\mathrm{ols}}}=CA^{({1}/{4})(1+z^2)}(A+D)^{({1}/{4})(7-z^2)+({1}/{2})mq_i}.\label{h.balance}
\end{equation}
In this balanced case, the $\hat A_i$ is unique provided $m>\frac{4+p}{1-q_i}$.
\end{theorem}

\begin{remark}\label{re1} Note that $h_{i;{\mathrm{gls}}}(A)$ does not have a
closed-form expression in $A$. But this is not an issue since finding a
root of the corresponding likelihood equation remains simple in this
case because the derivative of $\log[h_i(A)]$ has a closed-form
expression. Just like the standard residual likelihood, our adjusted
residual likelihood function could have multiple maxima in the general
balanced case. We refer to
Searle, Casella and McCulloch [(\citeyear{SCM92}), Section~8.1] who suggested
a way to search for the global maximum. In this connection, we refer to
\citet{GJ99} who proposed a method for testing for the global
maximum. Moreover, in order to reduce the number of iterations, we
suggest to use the simple ANOVA estimator of $A$ proposed by \citet{PR.1990} as an initial value.
\end{remark}

\begin{remark}\label{re2}
In a real data analysis, one should check
the condition $m>(4+p)/(1-h_i)$ for the existence of strictly positive
estimates $\hat{A}_{i;{\mathrm{gls}}}$ and $\hat{A}_{i;{\mathrm{gls}}}$. Under the
regularity conditions R2 and R3, the condition $m>(4+p)/(1-h_i)$
reduces to $m>m_0$, where $m_0$ is a fixed constant depending on $p$
and the leverages $q_i$. Thus, for sufficiently large $m$, this
condition does not pose any problem.
\end{remark}

\begin{remark}\label{re3} One might be tempted to treat our
adjustment factor $h_{i}(A)$ as a prior and conduct a regular
hierarchical Bayesian analysis. But $h_i(A)$ may not be treated as a
prior since in certain cases this leads
to an improper posterior distribution of $A$. To illustrate our point,
we consider the simple case: $D_i=D$ and
$h_{i}(A)=h_{i;{\mathrm{gls}}}(A)=h_{i;{\mathrm{ols}}}(A)$, $i=1,\ldots,m$. Since
\begin{eqnarray*}
h_{i}(A)L_{\mathrm{RE}}(A)&=&A^{(1+z^2)/4}(A+D)^{(7-z^2)/4+mq_i-m/2-p/2}
\\
&&{}\times\exp \biggl[-\frac{y^{\prime}(I-X^{\prime}(X^{\prime
}X)^{-1}X)y}{2(A+D)} \biggr]\bigl|X^{\prime}X\bigr|^{-1/2}C
\\
&\geq& 0,
\end{eqnarray*}
under the regularity conditions, and $\exp [-\frac{y^{\prime
}(I-X^{\prime}(X^{\prime}X)^{-1}X)y}{2(A+D)} ]$ and $A/(A+D)$ are
increasing monotone functions of $A$, there exists $s<\infty$ such that
\[
1- \exp \biggl[-\frac{y^{\prime}(I-X^{\prime}(X^{\prime
}X)^{-1}X)y}{2(s+D)} \biggr]<\frac{1}{2}
\]
and
\[
1-\frac{s}{s+D}<\frac{1}{2}.
\]
Using the above results, we have
%
\begin{equation}
\int_{0}^{\infty}h_{i}(A)L_{\mathrm{RE}}\,dA
\geq C\int_{s}^{\infty
}(A+D)^{2+1/2[mq_i+p]-m/2}\,dA,
\end{equation}
if $m>\frac{4+p}{1-q_i}$. Hence, if $-1\leq2+1/2[mq_i+p]-m/2\leq0$,
the right-hand side of the above equation is infinite, even if $m>\frac{4+p}{1-q_i}$.
Thus, in this case $h_i(A)$ cannot be treated as a prior since $\int_{0}^{\infty}h_{i}(A)L_{\mathrm{RE}}\,dA=\infty$ in case $-1\leq
2+1/2[mq_i+p]-m/2\leq0$.
\end{remark}
We now propose two empirical Bayes confidence intervals for $\theta_i$:
\[
I_i^{YL}(\hat A_{i;h})\dvtx \hat{
\theta}_i^{\mathrm{EB}}(\hat A_{i;h})\pm z_{\alpha/2}
\sigma_i(\hat A_{i;h}),
\]
where $h={\mathrm{gls}},{\mathrm{ols}}$. Since $\sigma_i(\hat A_{i;h})<
\sqrt{D_i}\ (h={\mathrm{gls}},{\mathrm{ols}})$, the length of our proposed Cox-type EB intervals, like
the original Cox EB interval $I_i^{\mathrm{Cox}}( \hat A_{\mathrm{ANOVA}})$, are always
shorter than that of the
direct interval $I_i^{D}$. The following theorem compares the lengths
of Cox EB confidence intervals of $\theta_i$ when $A$ is estimated by
$\hat A_{\mathrm{RE}}, \hat A_{i;{\mathrm{gls}}}$ and $\hat A_{i;{\mathrm{ols}}}$.

\begin{theorem}\label{th3}
Under the regularity conditions \textup{R2--R4} and $m>(4+p)/\break (1-q_i)$, we have
\[
\mbox{Length of } I_i^{\mathrm{Cox}}(\hat A_{\mathrm{RE}}) \le
\mbox{Length of } I_i^{YL}(\hat A_{i;{\mathrm{gls}}}) \le
\mbox{Length of } I_i^{YL}(\hat A_{i;{\mathrm{ols}}}).
\]
\end{theorem}

The following theorem provides the higher order asymptotic properties
of a general class of adjusted residual maximum likelihood estimators
of $A$.

\begin{theorem}\label{th4}
Under regularity conditions \textup{R1--R5}, we have:
\begin{longlist}[(ii)]
\item[(i)]$E[\hat{A}_{h_i}-A]=\frac{2}{\operatorname{tr}(V^{-1})}\tilde
{l}_{i,\mathrm{ad}}^{(1)}(A)+O(m^{-3/2})$,
\item[(ii)]$E(\hat{A}_{h_i}-A)^2=\frac{2}{\operatorname{tr}(V^{-1})}+O(m^{-3/2})$.
\end{longlist}
\end{theorem}

\begin{CT*}
Under regularity conditions \textup{R2--R5}, we have:
\begin{longlist}[(iii)]
\item[(i)] Both $\hat{A}_{i;{\mathrm{gls}}}$ and $\hat{A}_{i;{\mathrm{gls}}}$ are strictly
positive if $m>\frac{4+p}{1-q_i}$,
\item[(ii)]$E[\hat{A}_{i;{\mathrm{gls}}}-A]=\frac{2}{\operatorname{tr}(V^{-2})}\tilde
{l}_{i,\mathrm{ad};{\mathrm{gls}}}^{(1)}(A)+O(m^{-3/2})$,
\item[(iii)]$E[\hat{A}_{i;{\mathrm{ols}}}-A]=\frac{2}{\operatorname{tr}(V^{-2})}\tilde
{l}_{i,\mathrm{ad};{\mathrm{ols}}}^{(1)}(A)+O(m^{-3/2})$,
\item[(iv)]$E(\hat{A}_{i;h}-A)^2=\frac{2}{\operatorname{tr}(V^{-2})}+O(m^{-3/2})$,
\end{longlist}
where
\begin{eqnarray*}
\tilde{l}^{(1)}_{i;\mathrm{ad},{\mathrm{gls}}}&=&\frac{2}{A+D_i}+\frac
{(1+z^2)D_i}{4A(A+D_i)}+
\frac{1}{2}\operatorname{tr}\bigl(V^{-2}\bigr)x_{i}^{\prime}
\bigl(X^{\prime
}V^{-1}X\bigr)^{-1}x_i,
\\
\tilde{l}^{(1)}_{i;\mathrm{ad},{\mathrm{ols}}}&=&\frac{2}{A+D_i}+\frac
{(1+z^2)D_i}{4A(A+D_i)}+
\frac{1}{2}\operatorname{tr}\bigl(V^{-2}\bigr)x_{i}^{\prime}
\bigl(X^{\prime
}X\bigr)^{-1}X^{\prime}VX
\bigl(X^{\prime}X\bigr)^{-1}x_i.
\end{eqnarray*}
\end{CT*}

\begin{remark}\label{re4}
We reiterate that our true model variance
is $A$, which is not area specific (i.e., it does not depend on $i$).
However, unlike other likelihood based estimators of $A$, our theory
driven proposed adjusted REML estimators $\hat A_{i;{\mathrm{ols}}}$ and $\hat
A_{i;{\mathrm{gls}}}$ of $A$ are area and confidence level specific. We would like
to cite a similar situation that arises in the Bayesian small area
inference. For the same two level model, flat prior distribution on $A$
is widely accepted [see Morris and Tang (\citeyear{MT11})]. However, in order to
match the posterior variance with the classical MSE of EB with REML up
to the order $O(m^{-1})$, Datta, Rao and Smith (\citeyear{DRS.2005}) proposed a noncustomary
prior for $A$ that is area specific.
\end{remark}

\begin{remark}\label{re5} The area and confidence level specific
nature of our proposed estimators of a global parameter $A$ naturally
raises a concern that such proposed estimators may perform poorly when
compared to rival estimators of $A$. To address this issue,\vadjust{\goodbreak} first note
that the consistency of the new adjusted REML estimators $\hat
A_{i;{\mathrm{ols}}}$ and $\hat A_{i;{\mathrm{gls}}}$ of $A$ follows from part
(iv) of
the Corollary to Theorem~\ref{th4}. This is due to the fact that the leading term
in the right-hand side tends to 0 as $m$ tends to $\infty$, under the
regularity conditions R2--R5. This result also implies that MSEs of the
proposed estimators of $A$ are identical, up to the order $O(m^{-1})$,
to those of different likelihood based estimators of $A$ such as REML,
ML, different adjusted profile and residual maximum likelihood
estimators of \citet{L.L.2010} and \citet{YL}.
Moreover, while such an area and confidence level specific adjustment
causes the resulting proposed adjusted REML estimators to have more
bias than that of REML, the biases remain negligible and are of order
$O(m^{-1})$, same as the order of the bias of profile maximum
likelihood or adjusted profile maximum likelihood estimators of $A$
proposed by \citet{L.L.2010} and \citet{YL}.
Basically, we introduce this slight bias in $\hat A_{i;{\mathrm{ols}}}$ and $\hat
A_{i;{\mathrm{gls}}}$ in order to achieve the desired low coverage error property
while maintaining length always shorter than that of the corresponding
direct confidence interval.
\end{remark}

\begin{remark}\label{re6} Using the Corollary to Theorem~\ref{th4} and the
mathematical tools used in \citet{L.L.2010}, we obtain the
following second-order approximation to the mean squared error (MSE) of
$\hat{\theta}_i^{\mathrm{EB}}(\hat A_{i;{\mathrm{gls}}})$:
\[
\operatorname{MSE}\bigl[\hat{\theta}_i^{\mathrm{EB}}(\hat
A_{i;{\mathrm{gls}}})\bigr] =g_{1i}(A) + g_{2i}(A) +
g_{3i}(A)+o\bigl(m^{-1}\bigr),
\]
where
$g_{1i}(A){=}\frac{AD_i}{A+D_i}, g_{2i}(A){=}\frac{ D_i^{2} }{ (A+D_i)^2
}\operatorname{ Var}(x_i'\hat\beta){=}\frac{ D_i^{2} }{ (A+D_i)^2 }x_i'
(\sum_{j=1}^m \frac{x_jx_j' }{A+D_j})^{-1}\times x_i$, and
$g_{3i} (A)=\frac{2D_i^2}{(A+D_i)^3} \{\sum_{j=1}^m \frac
{1}{(A+D_j)^2}
\}^{-1}$.
Thus, in terms of MSE criterion, $\hat{\theta}_i^{\mathrm{EB}}[\hat A_{i;{\mathrm{gls}}}]$
is equally efficient, up to the order $O(m^{-1})$, as the empirical
Bayes estimators of $\theta_i$ that use standard REML, PML and the
adjusted PML and REML estimators of $A$ proposed by \citet{L.L.2010} and \citet{YL}.

We note that
\[
\operatorname{MSE}\bigl[\hat{\theta}_i^{\mathrm{EB}}(\hat
A_{i;{\mathrm{ols}}})\bigr] =g_{1i}(A) + g_{2i;{\mathrm{ols}}}(A) +
g_{3i}(A)+o\bigl(m^{-1}\bigr),
\]
where
$g_{2i;{\mathrm{ols}}}(A) =\frac{ D_i^{2} }{ (A+D_i)^2 }x_i^{\prime}(X^{\prime
}X)^{-1}X^{\prime}VX(X^{\prime}X)^{-1}x_i\ge\frac{ D_i^{2} }{
(A+D_i)^2 }\times\break  x_i^{\prime}(X^{\prime}V^{-1}X)^{-1}x_i$. Thus, in terms of
higher order asymptotics $\hat{\theta}_i^{\mathrm{EB}}(\hat A_{i;{\mathrm{ols}}})$ is less
efficient than $\hat{\theta}_i^{\mathrm{EB}}(\hat A_{i;{\mathrm{gls}}})$.
\end{remark}

\begin{remark}\label{re7} We suggest the following second-order
unbiased estimator of $\operatorname{MSE}[\hat{\theta}_i^{\mathrm{EB}}(\hat A_{i;{\mathrm{gls}}})]$:
\[
\mathrm{mse}_i=g_{1i}(\hat{A}_{i;{\mathrm{gls}}})+g_{2i}(
\hat {A}_{i;{\mathrm{gls}}})+2g_{3i}(\hat{A}_{i;{\mathrm{gls}}})-\bigl[
\hat{B}_{i}(\hat {A}_{i;{\mathrm{gls}}})\bigr]^2\widehat{
\operatorname{Bias}}(\hat{A}_{i;{\mathrm{gls}}}),
\]
where $\hat{B}_{i}(\hat{A}_{i;{\mathrm{gls}}})=\frac{D_i}{D_i+\hat{A}_{i;{\mathrm{gls}}}}$,
and $\widehat{\operatorname{Bias}}(\hat{A}_{i;{\mathrm{gls}}})=\frac{2}{\operatorname{tr}(V^{-2})}\tilde
{l}_{i,\mathrm{ad};{\mathrm{gls}}}^{(1)}(\hat{A}_{i;{\mathrm{gls}}})$.
We provide expressions for the second-order MSE approximation and the
second-order unbiased estimator of $\operatorname{MSE}[\hat{\theta
}_i^{\mathrm{EB}}(\hat A_{i;{\mathrm{gls}}})]$ for the benefit of researchers interested in
such expressions. However, for the purpose of point estimation and the
associated second-order unbiased MSE estimators, we recommend the
estimators proposed by \citet{YL}. We recommend the
use of $\hat{A}_{i;{\mathrm{gls}}}$ only for the construction of second-order
efficient Cox-type EB confidence intervals.
\end{remark}

\section{A Monte Carlo simulation study}\label{sec4}

In this section, we design a Monte Carlo simulation study to
compare finite sample performances of the following confidence
intervals of $\theta_i$ for the Fay--Herriot model: direct, Cox-type
EB using (i) REML estimator of $A$ (Cox.RE), (ii) estimator of $A$
proposed by \citet{WF.2003} (Cox.WF), (iii) estimator of $A$
proposed by \citet{L.L.2010} (Cox.LL), parametric bootstrap EB
confidence interval of Chatterjee, Lahiri and Li (\citeyear{CLL.2008}) using Li--Lahiri
estimator of $A$ (CLL.LL), our proposed Cox-type EB confidence
intervals using GLS estimator of $\beta$ (Cox.YL.GLS) and OLS estimator
of $\beta$ (Cox.YL.OLS). In Section~\ref{sec4.1}, we consider a Fay--Herriot
model with a common mean as in Datta, Rao and Smith (\citeyear{DRS.2005}) and
Chatterjee, Lahiri and Li (\citeyear{CLL.2008}). In Section~\ref{sec4.2}, we consider a Fay--Herriot model with
one auxiliary variable in order to examine the effect of different
leverage and sampling variance combinations on the coverage and average
length of different confidence intervals of a small area mean.

\subsection{The Fay--Herriot model with a common mean}\label{sec4.1}
Throughout this subsection, we assume a common mean $x_i'\beta=0$,
which is estimated using data as in other papers on small area
estimation. Specifically, we generate $R=10^4$ independent replicates
$\{y_{i}, v_i, i=1,\ldots,m\}$ using the following Fay--Herriot model:
\[
y_{i}=v_{i}+e_{i},
\]
where $v_{i}$ and $e_{i}$ are mutually independent with $v_{i}\stackrel
{\mathrm{i.i.d.}}{\sim}N(0,A)$, $e_{i}\stackrel{\mathrm{ind}}{\sim}N(0,D_i),\break  i=1,\ldots,m$.
We set $A=1$. For the parametric bootstrap method, we consider $B=6000$
bootstrap samples.

In the unbalanced case, for $m=15$, we consider five groups,
say $G\equiv(G_{1}, G_{2}, G_{3}, G_{4}, G_{5})$, of small areas, each
with three small areas, such that the sampling variances $D_i$ are the
same within a given area. We consider the following two patterns of the
sampling variances: (a) $(0.7,0.6,0.5,0.4,0.3)$ and (b)~$(4.0, 0.6,
0.5, 0.4, 0.1)$. Note that in pattern (a) all areas have sampling
variances less than $A$. In contrast, in pattern (b), sampling
variances of all but one area are less than $A$. The patterns (a) and
(b) correspond to the sampling variance patterns (a) and (c) of Datta, Rao and Smith (\citeyear{DRS.2005}).

The simulation results are displayed in Table~\ref{tab1}. First note that while
the direct method attains the nominal coverage most of the time it has
the highest length compared to the other methods considered. The
interval Cox.RE cuts down the length of the direct method considerably
at the expense of undercoverage, which is more severe for pattern (b)
than pattern (a). This could be due to the presence of three outlying
areas (i.e., with respect to the sampling variances) in $G_1$. The
intervals Cox.WF and Cox.LL improve on Cox.RE as both use strictly
positive consistent estimators of $A$.
Our new methods---Cox.YL.GLS and Cox.YL.OLS---and CLL.LL perform
very well in terms of coverage although CLL.LL is showing a slight
undercoverage.
The CLL.LL method is slightly better than ours in terms of average
length although we notice that in some simulation replications the
length of the parametric bootstrap EB confidence interval is larger
than that of the direct.\looseness=-1

\subsection{Effect of leverage and sampling variance in a Fay--Herriot
model with one auxiliary variable}\label{sec4.2}
We generate $R=10^4$ independent replicates $\{y_{i}, v_i, i=1,\ldots,m\}$
using the following Fay--Herriot model:
\[
y_{i}=x_i\beta+v_{i}+e_{i},\vadjust{\goodbreak}
\]
where $v_{i}$ and $e_{i}$ are mutually independent with $v_{i}\stackrel
{\mathrm{i.i.d.}}{\sim}N(0,A)$, $e_{i}\stackrel{\mathrm{ind}}{\sim}N(0,D_i),\break  i=1,\ldots,m$.
We set $A=1$. For the parametric bootstrap method, we consider $B=6000$
bootstrap samples.

In this subsection, we examine the effects of leverage and sampling
variance on different confidence intervals for $\theta_i$. We consider
six different (leverage, sampling variance) patterns of the first area
using leverages $(0.07,0.22,\break 0.39)$ and sampling variances
$D_1=(1,5,10)$. For the remaining 14 areas, we assume equal small
sampling variances $D_j=0.01, j\ge2$ and same leverage. Since the total
leverage for all the areas must be 1, we obtain the common leverage for
the other areas from the knowledge of leverage for the first area.

In Table~\ref{tab2}, we report the coverages and average lengths for all the
competing methods for the first area for all the six patterns. We do
not report the results for the remaining 14 areas since they are
similar, as expected, due to small sampling variances in those areas.
The use of strictly positive consistent estimators of $A$ such as WF
and LL help bringing coverage of the Cox-type EB confidence interval
closer to the nominal coverage of $95\%$ than the one based on REML.
For large sampling variances and leverages, the Cox-type EB confidence
intervals based on REML, WF and LL methods have generally shorter
length than ours or parametric bootstrap confidence interval but only
at the expense of severe undercoverage. Our simulation results show
that our proposed Cox.YL.GLS could perform better than Cox.YL.OLS and
is very competitive to the more computer intensive CLL.LL method.

\begin{sidewaystable}
\tablewidth=\textwidth
\caption{Simulation results for Section~\protect\ref{sec4.1}: Simulated coverage and
average length (in parenthesis) of different confidence intervals of
small area means; nominal coverage is $95\%$}
\label{tab1}
\begin{tabular*}{\textwidth}{@{\extracolsep{\fill}}lccccccccccccccc@{}}
\hline
\textbf{Pattern} & \textbf{G} & \multicolumn{2}{c}{\textbf{Cox.WF}} & \multicolumn{2}{c}{\textbf{Cox.RE}}
& \multicolumn{2}{c}{\textbf{Cox.LL}} & \multicolumn{2}{c}{\textbf{CLL.LL}} &
\multicolumn{2}{c}{\textbf{Cox.YL.GLS}} & \multicolumn{2}{c}{\textbf{Cox.YL.OLS}} &
\multicolumn{2}{c@{}}{\textbf{Direct}} \\
\hline
a & 1 & 90.6 & (2.4) & 90.4 & (2.4) & 94.2 & (2.6) & 94.9 & (2.7) &
95.3 & (2.8) & 95.3 & (2.8) & 95.1 & (3.3) \\
& 2 & 91.2 & (2.3) & 90.8 & (2.3) & 94.3 & (2.5) & 94.9 & (2.5) & 95.3
& (2.6) & 95.3 & (2.6) & 94.9 & (3.0) \\
& 3 & 91.5 & (2.1) & 90.8 & (2.1) & 94.2 & (2.3) & 94.9 & (2.4) & 95.3
& (2.4) & 95.3 & (2.4) & 95.1 & (2.8) \\
& 4 & 91.8 & (2.0) & 91.2 & (2.0) & 94.3 & (2.1) & 94.9 & (2.2) & 95.2
& (2.2) & 95.3 & (2.2) & 95.2 & (2.5) \\
& 5 & 92.4 & (1.8) & 92.1 & (1.8) & 94.7 & (1.9) & 95.1 & (1.9) & 95.5
& (2.0) & 95.5 & (2.0) & 95.1 & (2.1) \\[3pt]
b & 1 & 88.3 & (3.3) & 88.1 & (3.3) & 93.7 & (3.8) & 94.6 & (4.0) &
95.6 & (4.3) & 95.9 & (4.3) & 94.8 & (7.8) \\
& 2 & 90.9 & (2.3) & 90.0 & (2.3) & 93.9 & (2.5) & 94.7 & (2.5) & 95.2
& (2.6) & 95.3 & (2.6) & 94.9 & (3.0) \\
& 3 & 91.2 & (2.1) & 90.2 & (2.1) & 93.9 & (2.3) & 94.7 & (2.4) & 95.0
& (2.5) & 95.2 & (2.5) & 95.1 & (2.8) \\
& 4 & 91.7 & (2.0) & 90.9 & (2.0) & 94.4 & (2.1) & 95.0 & (2.2) & 95.3
& (2.2) & 95.4 & (2.3) & 95.0 & (2.5) \\
& 5 & 93.8 & (1.1) & 93.1 & (1.1) & 94.8 & (1.2) & 94.9 & (1.2) & 95.0
& (1.2) & 95.0 & (1.2) & 94.9 & (1.2) \\
\hline
\end{tabular*}\vspace*{9pt}
%
%
%
\tablewidth=\textwidth
\caption{Simulation results for Section~\protect\ref{sec4.2}: Simulated coverage and
average length (in parenthesis) of different confidence intervals for
the first small area mean for different combinations of leverage and
sampling variance of the first area; nominal coverage is $95\%$}
\label{tab2}
\begin{tabular*}{\textwidth}{@{\extracolsep{\fill}}lccccccccccccccc@{}}
\hline
\textbf{Leverage} & $\bolds{D_1}$ & \multicolumn{2}{c}{\textbf{Cox.WF}} & \multicolumn
{2}{c}{\textbf{Cox.RE}} & \multicolumn{2}{c}{\textbf{Cox.LL}} & \multicolumn
{2}{c}{\textbf{CLL.LL}} & \multicolumn{2}{c}{\textbf{Cox.YL.{{gls}}}} & \multicolumn
{2}{c}{\textbf{Cox.YL.{{ols}}}} & \multicolumn{2}{c}{\textbf{Direct}} \\
\hline
0.39 & 10 & 78.1 & (3.2) & 85.3 & (3.6) & 88.0 & (3.9) & 94.7 & (5.0) &
98.0 & (6.9) & 98.3 & (8.1) & 95.1 & (12.4) \\
& \phantom{0}5 & 81.5 & (3.2) & 86.6 & (3.5) & 89.1 & (3.7) & 95.0 & (4.6) & 97.0
& (5.8) & 97.3 & (6.2) & 94.9 & \phantom{0}(8.8) \\
& \phantom{0}1 & 89.7 & (2.7) & 90.0 & (2.7) & 91.3 & (2.8) & 94.9 & (3.2) & 95.3
& (3.4) & 95.4 & (3.4) & 94.8 & \phantom{0}(3.9) \\[3pt]
0.22 & 10 & 84.0 & (3.4) & 89.7 & (3.7) & 92.2 & (3.9) & 95.3 & (4.5) &
96.7 & (5.0) & 98.5 & (5.7) & 94.9 & (12.4) \\
& \phantom{0}5 & 85.8 & (3.3) & 89.9 & (3.5) & 91.9 & (3.8) & 95.0 & (4.2) & 96.0
& (4.6) & 97.1 & (4.9) & 95.0 &\phantom{0}(8.8) \\
& \phantom{0}1 & 91.6 & (2.7) & 91.9 & (2.7) & 93.0 & (2.8) & 94.9 & (3.0) & 95.5
& (3.2) & 95.5 & (3.2) & 95.2 & \phantom{0}(3.9) \\[3pt]
0.07 & 10 & 87.2 & (3.5) & 92.2 & (3.7) & 94.2 & (3.9) & 95.3 & (4.1) &
95.7 & (4.2) & 96.1 & (4.3) & 95.0 & (12.4) \\
& \phantom{0}5 & 89.2 & (3.4) & 92.7 & (3.5) & 94.4 & (3.7) & 95.5 & (3.9) & 95.8
& (4.0) & 95.9 & (4.0) & 94.8 & \phantom{0}(8.8) \\
& \phantom{0}1 & 93.3 & (2.7) & 93.3 & (2.7) & 94.4 & (2.8) & 95.2 & (2.9) & 95.4
& (3.0) & 95.4 & (3.0) & 95.2 & \phantom{0}(3.9) \\
\hline
\end{tabular*}
\end{sidewaystable}

\section{Concluding remarks}\label{sec5}
In this paper, we put forward a new simple approach for constructing
second-order efficient empirical Bayes confidence interval for a small
area mean using a carefully devised adjusted residual maximum
likelihood estimator of the model variance in the well-known Cox
empirical Bayes confidence interval. Our simulation results show that
the proposed method performs much better than the direct or Cox EB
confidence intervals with different standard likelihood based
estimators of the model variance. In our simulation, the parametric
bootstrap empirical Bayes confidence interval also performs well and it
generally produces intervals shorter than direct confidence intervals
on the average. However, to the best of our knowledge, there is no
analytical result that shows that the parametric bootstrap empirical
Bayes confidence interval is always shorter than the direct interval.
In fact, in our simulation we found cases where the length of
parametric bootstrap empirical Bayes confidence interval is higher than
that of the direct. In order to obtain good parametric bootstrap
empirical Bayes confidence intervals, choices of the estimator of $A$
and the bootstrap replication $B$ appear to be important. To limit the
computing time, we have considered a simple simulation setting with
$m=15$. During the course of our investigation, we feel the need for
developing an efficient computer program that allows evaluation of
parametric bootstrap empirical Bayes confidence intervals in a large
scale simulation environment. Until the issues raised for the
parametric bootstrap empirical Bayes confidence interval method are
resolved, our proposed simple second-order efficient empirical Bayes
confidence interval could serve as a promising method. The results
presented in this paper is for the well-known Fay--Herriot model. It is
not clear at this time how the results will extend to a general class
of small area models---this will be a challenging topic for future research.

\begin{appendix}\label{app}
\section{}

In this appendix, we provide an outline of proofs of different
results presented in the paper. In order to facilitate the review, we
supply a detailed proof of Theorem~\ref{th4} in the supplementary material
[\citet{supp}].

\begin{pf*}{Proof of Theorem~\ref{th1}}
For notational simplicity, we set $\hat A_{h_i}\equiv\hat A$
throughout the \hyperref[app]{Appendix}. Define
\[
G_i(z,y)= z[\hat{\sigma}_i/\sigma_i-1]+
\bigl\{(B_i-\hat{B}_i) \bigl(y_i-x^{\prime
}_i
\beta\bigr)+\hat{B}_i \bigl[x^{\prime}_i(\hat{
\beta}-\beta)\bigr]\bigr\}/\sigma_i.
\]
Using calculations similar to the ones
Chatterjee, Lahiri and Li (\citeyear{CLL.2008}), we have
%
\begin{eqnarray}\label{app.1}
&&P\bigl[\theta_{i}\leq\hat{\theta}^{\mathrm{EB}}_i(
\hat{A})+z\hat{\sigma}_i\bigr]
\nonumber
\\
&&\qquad=\Phi(z)+\phi(z)E \biggl[G_i(z,y) -\frac{z}{2}G_i^{2}(z,y)
\biggr]\\
&&\qquad\quad{}+ \frac{1}{2}E \biggl[\int_{z}^{z+G_i(z,y)}
\bigl(z+G_i(z,y)-t\bigr)^2 \bigl(t^2-1\bigr)
\phi(t) \,dt \biggr].
\nonumber
\end{eqnarray}
We shall first show that the third term of the right-hand side of (\ref
{app.1}) is of order $O(m^{-3/2})$. To this end,
using
\[
0\leq\bigl|z+G_i(z,y)-t\bigr|\leq\bigl|G_i(z,y)\bigr| \quad{\mbox{and}}\quad
\bigl(t^2-1\bigr)\phi(t)\leq 2\phi(\sqrt{3}),
\]
in $t\in(z,z+G_i(z,y))$, we have
\begin{eqnarray*}
\mathrm{third\ term\ of\ (\ref{app.1})}&\leq&\frac{1}{2}E \biggl[\int
_{z}^{z+G_i(z,y)}\bigl(z+G_i(z,y)-t
\bigr)^2 \bigl|\bigl(t^2-1\bigr)\phi(t)\bigr| \,dt \biggr]
\\
&\leq&C\phi(\sqrt{3})E\bigl[G_i^{3}(z,y)\bigr].
\end{eqnarray*}

Setting $\sigma_i^2=S_i$ and using the Taylor series expansion, we have
\[
\hat{\sigma}_i(\hat{S}_i)- \sigma_i(S_i)=
\tfrac{1}{2}S_i^{- {1}/{2}}(\hat{S}_i-S_i)
-\tfrac{1}{8}S_i^{-{3}/{2} }(\hat{S}_i-S_i)^2+O_p
\bigl(|\hat{S}_i-S_i|^{3}\bigr),
\]
so that
\[
\frac{\hat{\sigma}_i(\hat{S}_i)}{\sigma_i(S_i)}-1=\frac{1}{2S_i} (\hat{S}_i-S_i)
-\frac{1}{8S_i^2}(\hat{S}_i-S_i)^2+R_{A1}.
\]

Using
\begin{eqnarray*}
\hat{B}_i-B_i&=&-(\hat{A}-A)\frac{D_i}{(A+D_i)^2}+(
\hat{A}-A)^2 \frac
{D_i}{(A+D_i)^3}+R_{A2},
\\
\hat{\sigma}_i^2-\sigma_i^2&=&(
\hat{A}-A)\frac{D_i^2}{(A+D_i)^2}-(\hat {A}-A)^2 \frac{D_i^2}{(A+D_i)^3}+R_{A3},
\end{eqnarray*}
we can write $G_i(z,y)=G_{1i}(y)+G_{2i}(z,y)$, where
\begin{eqnarray*}
G_{1i}(y)&=&\frac{1}{\sqrt{m}}\hat{u}_{1i}+\frac{1}{m}
\hat{u}_{2i}+R_{A4},
\\
G_{2i}(z,y)&=&z \biggl[\frac{1}{\sqrt{m}}\hat{v}_{1i}+
\frac{1}{m}\hat {v}_{2i} \biggr] +R_{A5},
\end{eqnarray*}
with
\begin{eqnarray*}
\hat{u}_{1i}&=&\sqrt{m} \sigma_i^{-1}
\biggl[B_ix_i^{\prime}(\hat{\beta }-\beta)+(
\hat{A}-A)\frac{D_i}{(A+D_i)^2}\bigl(y_i-x_i^{\prime}
\beta\bigr) \biggr],
\\
\hat{u}_{2i}&=&m \sigma_i^{-1} \biggl[-(
\hat{A}-A)^2\frac
{D_i}{(A+D_i)^3}\bigl(y_i-x_i^{\prime}
\beta\bigr)\\
&&\hspace*{32pt}{}+(\hat{A}-A)\frac
{D_i}{(A+D_i)^2}B_ix_i^{\prime}(
\hat{\beta}-\beta) \biggr],
\\
\hat{v}_{1i}&=&\sqrt{m} \frac{B_i^2}{2\sigma_i^2}(\hat{A}-A),
\\
\hat{v}_{2i}&=&m \biggl[-\frac{1}{2\sigma_i^2}\frac{B_i^2}{A+D_i}(\hat
{A}-A)^2-\frac{1}{8\sigma_i^4}(\hat{A}-A)^2B_i^4
\biggr].
\end{eqnarray*}

Using the fact that $E[|\hat{A}-A|^k]=O(m^{-3/2})$ for $k\geq3$ [this
can be proved using the mathematical tools used in \citet{L.L.2010}
and Das, Jiang and Rao (\citeyear{Das})], we have, for $k=1,2,3,4,5$ and large $m$,
\begin{eqnarray*}
E\bigl[|R_{Ak}|\bigr]&\leq& C E \bigl[|\hat{A}-A|^3 \bigr]=O
\bigl(m^{-3/2}\bigr),
\\
\bigl|\mathrm{third\ term\ of\ (\ref{app.1})}\bigr|&\leq& C\phi(\sqrt {3})E
\bigl[\bigl|G_i^{3}(z,y)\bigr|\bigr]\leq C E \bigl[|
\hat{A}-A|^3 \bigr]=O\bigl(m^{-3/2}\bigr),
\end{eqnarray*}
where $C$ is a generic constant.

We also note that
\[
E\bigl[G_i(z,y)\bigr]=m^{-1/2}E[\hat{u}_{1i}+z
\hat{v}_{1i}]+m^{-1}E[\hat {u}_{2i}+z
\hat{v}_{2i}]+O\bigl(m^{-3/2}\bigr),
\]
so that
\[
\mbox{the right-hand side of } (\ref{app.1})=\Phi(z)+\phi(z)E
\biggl[G_i(z,y) -\frac{z}{2}G_i^{2}(z,y)
\biggr]+O\bigl(m^{-3/2}\bigr).
\]

Similarly,
%
\begin{eqnarray}
P\bigl(\hat{\theta}_{i}^{\mathrm{EB}}-z\hat{\sigma}_i
\leq\theta_i\bigr)&=&\Phi(-z)+\phi (-z)E \biggl[G_i(-z,y)
+\frac{z}{2}G_i^{2}(-z,y) \biggr]
\nonumber
\\[-8pt]
\\[-8pt]
\nonumber
&&{}+ O
\bigl(m^{-3/2}\bigr),
\end{eqnarray}
so that using
\begin{eqnarray*}
&&G_i(z,y)-G_i(-z,y)-\frac{z}{2}
\bigl[G_i^{2}(z,y)+G_i^{2}(-z,y)
\bigr]
\\
&&\qquad=2G_{2i}(z,y)-\frac{z}{2}\bigl[G_{1i}^2(y)+G_{2i}^2(z,y)
\bigr]
\\
&&\qquad=\frac{2z}{\sqrt{m}}{\hat{v}_{1i}}+\frac{z}{m}\bigl\{2
\hat{v}_{2i}-\hat {u}_{1i}^2-z^2
\hat{v}_{1i}^2\bigr\}+R_{A6},
\end{eqnarray*}
where $E[|R_{A6}|]=O(m^{-3/2})$ since $E[|\hat{A}-A|^k]=O(m^{-3/2})
\mbox{ for $k\ge3$}$.

We have
\begin{eqnarray*}
&&P \bigl\{\theta_i \in I_{cox}(\hat{A}) \bigr\}
\\
&&\qquad=\Phi(z)-\Phi(-z)+\phi(z)E\bigl[G_i(z,y)-G_i(-z,y)
\bigr]\\
&&\qquad\quad{}-\frac{z}{2}\phi (z)E\bigl[G_i^{2}(z,y)+G_i^{2}(-z,y)
\bigr]+O\bigl(m^{-3/2}\bigr)
\\
&&\qquad=1-\alpha+z\phi(z) \bigl\{m^{-{1}/{2}}E[2\hat {v}_{1i}]+m^{-1}E
\bigl[2\hat{v}_{2i}-\hat{u}_{1i}^2-z^2
\hat{v}_{1i}^2\bigr] \bigr\}\\
&&\qquad\quad{}+O\bigl(m^{-3/2}\bigr).
\end{eqnarray*}
Using Lemma~\ref{le1}, given below, and considerable algebra, we show that
\[
a_i= E\bigl[2\hat{v}_{2i}-\hat{u}_{1i}^2-z^2
\hat{v}_{1i}^2\bigr] \quad\mbox{and}\quad b_i=2\sqrt{m}E[
\hat{v}_{1i}].
\]
This completes the proof of equation (\ref{PI.cox}).
\end{pf*}

\begin{lem}\label{le1}
Under the regularity conditions \textup{R1--R5}, we have
%
\begin{eqnarray}\qquad
E\bigl[\hat{v}_{1i}^2(\hat{A})\bigr]&=&\frac{m}{\operatorname{tr}(V^{-2})}
\frac
{D_i^2}{2A^2(A+D_i)^2}+O\bigl(m^{-1/2}\bigr),\label{v1.2}
\\[1pt]
\label
{v2}E\bigl[\hat{v}_{2i}(\hat{A})\bigr]&=&-\frac{m}{\operatorname{tr}(V^{-2})} \biggl[
\frac
{D_i}{A(A+D_i)^2}+\frac{D_i^2}{4A^2(A+D_i)^2} \biggr]
\nonumber
\\[-7pt]
\\[-7pt]
\nonumber
&&{}+O\bigl(m^{-1/2}
\bigr),
\\[1pt]
\label{u1.2}
E\bigl[\hat{u}_{1i}^{2}(\hat{A})\bigr]&=&m\frac{D_i}{A(A+D_i)}
\biggl[E\bigl[\bigl\{ x_i^{\prime}\tilde{\beta}- \beta)\bigr
\}^2\bigr]
+\frac{D_i}{A(A+D_i)^2}\frac
{2}{\operatorname{tr}(V^{-2})} \biggr]
\nonumber
\\[-7pt]
\\[-7pt]
\nonumber
&&{}+O
\bigl(m^{-1/2}\bigr),
\\[1pt]
E\bigl[\hat{v}_{1i}(\hat{A})\bigr]&=&\frac{\sqrt{m}}{\operatorname{tr}(V^{-2})}
\frac
{D_i}{A(A+D_i)}\tilde{l}_{i;\mathrm{ad}}^{(1)}+O\bigl(m^{-1}
\bigr).\label{v1.adRE}
\end{eqnarray}
\end{lem}

\begin{pf*}{Proof of Theorem~\ref{th2}}
First note that solution of $h_i(A)$ to the differential equation (\ref
{diffeqn}) depends on whether the OLS or GLS method is used to estimate
${\beta}$. Also note that the solution of $h_i(A)$ for the OLS case
does not follow as a special case of GLS. Thus, we treat these two
cases separately. The balanced case, that is, equation (\ref
{h.balance}) follows from (\ref{h.unbalance.gls}) or (\ref
{h.unbalance.ols}).

\textit{Case \textup{1:} Unbalanced case: OLS \textup{[}proof of equation \textup{(\ref
{h.unbalance.ols})]}}

From equation (\ref{diffeqn}), we have
\[
\tilde{l}_{i;\mathrm{ad}}^{(1)}(A)=\frac{2}{A+D_i}+
\frac
{(1+z^2)D_i}{4A(A+D_i)}+\frac{1}{2}x_i^{\prime}
\bigl(X^{\prime
}X\bigr)^{-1}X^{\prime}VX
\bigl(X^{\prime}X\bigr)^{-1}x_i\operatorname{tr}
\bigl(V^{-2}\bigr).
\]
Therefore,
\begin{eqnarray*}
\tilde{l}_{i;\mathrm{ad}}(A)&=&\int\tilde{l}_{i;\mathrm{ad}}^{(1)}\,dA
\\[1pt]
&=&2\log(A+D_i)+\frac{(1+z^2)}{4}
\log \biggl(\frac{A}{A+D_i} \biggr)\\[1pt]
&&{}+\frac
{1}{2}x_i^{\prime}
\bigl(X^{\prime}X\bigr)^{-1}X^{\prime}J X
\bigl(X^{\prime
}X\bigr)^{-1}x_i+C
\\[1pt]
&=&2\log(A+D_i)+\frac{(1+z^2)}{4}\log \biggl(\frac{A}{A+D_i}
\biggr)
\\[1pt]
&&{}+\frac{1}{2}x_i^{\prime}\bigl(X^{\prime}X
\bigr)^{-1}X^{\prime} \bigl[-V\operatorname{tr}\bigl(V^{-1}\bigr)+\operatorname{tr}
\bigl(V^{-1}\bigr)+C \bigr] X\bigl(X^{\prime}X
\bigr)^{-1}x_i+C
\\[1pt]
&=&2\log(A+D_i)+\frac{1+z^2}{4}\log \biggl(\frac{A}{A+D_i}
\biggr)
\\
&&{}-\frac{1}{2}x_i^{\prime}\bigl(X^{\prime}X
\bigr)^{-1}X^{\prime}V X\bigl(X^{\prime
}X
\bigr)^{-1}x_i\operatorname{tr}\bigl(V^{-1}\bigr)\\
&&{}+
\frac{1}{2}q_i \Biggl[\sum_{i=1}^m
\log (A+D_i) \Biggr]+C.
\end{eqnarray*}
In addition,
\[
J=\operatorname{diag} \biggl(\int(A+D_1 )\operatorname{tr}\bigl(V^{-2}\bigr)\,dA,\ldots,
\int(A+D_m)\operatorname{tr}\bigl(V^{-2}\bigr)\,dA \biggr).
\]
Equation (\ref{h.unbalance.ols}) follows noting that $h_{i}(A)=\exp
[\tilde{l}_{i;\mathrm{ad}}(A)]$.

\textit{Case \textup{2:} Unbalanced case: GLS \textup{[}proof of \textup{(\ref
{h.unbalance.gls})]}}

Solving equation (\ref{diffeqn}) for $\tilde{l}_{i;\mathrm{ad}}^{(1)}(A)$, we get
\[
\tilde{l}_{i;\mathrm{ad}}^{(1)}(A)=\frac{2}{A+D_i}+
\frac
{(1+z^2)D_i}{4A(A+D_i)}+\frac{1}{2}x_i^{\prime}
\bigl(X^{\prime}V^{-1}X\bigr)^{-1} x_i\operatorname{tr}
\bigl(V^{-2}\bigr).
\]
Thus,
\begin{eqnarray*}
\tilde{l}_{i;\mathrm{ad}}(A)&=&\int\tilde{l}_{i;\mathrm{ad}}^{(1)}\,dA
\\
&=&\int\frac{2}{A+D_i}\,dA+\int\frac{(1+z^2)D_i}{4A(A+D_i)}\,dA\\
&&{}+\frac
{1}{2}\int
x_i^{\prime}\bigl(X^{\prime}V^{-1}X
\bigr)^{-1} x_i\operatorname{tr}\bigl(V^{-2}\bigr) \,dA
\\
&=&2\log(A+D_i)+\frac{1}{4}\bigl(1+z^2
\bigr)D_i\log\frac{A}{A+D_i}+\frac
{1}{2}K+C,\qquad \mbox{say,}
\end{eqnarray*}
where $K=\int x_i^{\prime}(X^{\prime}V^{-1}X)^{-1} x_i\operatorname{tr}(V^{-2}) \,dA$.

We now prove part (ii) of the theorem. To this end, note that the
adjusted maximum residual likelihood estimator of $A$ with the
adjustment factor (\ref{h.balance}) is obtained as a solution of
\begin{eqnarray*}
&&l_{\mathrm{RE}}^{(1)}+\tilde{l}_{i,\mathrm{ad}}^{(1)}=0
\\
&&\qquad
\Longleftrightarrow f(A)\equiv\bigl\{-2(m-p)+8+2m q_i\bigr
\}A^2
\\
&&\hspace*{50pt}\qquad\quad{}+\bigl\{2y^{\prime}\bigl(I_m-X\bigl(X^{\prime}X
\bigr)^{-1}X^{\prime
}\bigr)y-2(m-p)D+8D\\
&&\hspace*{201pt}{}+\bigl(1+z^2
\bigr)D+2mD q_i\bigr\}A
\\
&&\hspace*{50pt}\qquad\quad{}+\bigl(1+z^2\bigr)D^2=0.
\end{eqnarray*}
Therefore, under strict positiveness of the solution and $m>\frac
{4+p}{1-q_i}$, $f(A)$ is a quadratic and concave function of $A$.
Thus, due to $f(0)>0$, there is a unique and strictly positive adjusted
residual maximum likelihood estimator of $A$ in the balanced case.
\end{pf*}

\begin{pf*}{Proof of Theorem~\ref{th3}}
Note that the length of the Cox-type EB confidence interval of $\theta
_i$ is given by $2\sigma(\hat A_i)$, where $\sigma(\hat A_i)=\sqrt
{\frac{\hat A_i D_i}{\hat A_i+D_i}}$ and $\hat A_i$ is an estimator of
$A$ used to construct an empirical Bayes confidence interval for $\theta
_i$. We show that among the three intervals considered the length of
the Cox EB confidence interval is the shortest when $\hat A_{\mathrm{RE}}$ is
used to estimate $A$, followed by $\hat A_{i,{\mathrm{gls}}}$, and $\hat
A_{i,{\mathrm{ols}}}$. Since $\sigma(\hat A_i)$ is a monotonically increasing
function of $\hat A_i$, it suffices to show that
\[
\hat{A}_{\mathrm{RE}}\leq\hat{A}_{i,{\mathrm{gls}}}\leq\hat{A}_{i,{\mathrm{ols}}}.
\]

Note that
\begin{eqnarray*}
l_{\mathrm{RE}}^{(1)}(\tilde{A}_{\mathrm{RE}})&=&0,
\\
l_{\mathrm{RE}}^{(1)}(\hat{A}_{i,{\mathrm{gls}}})+\tilde{l}^{(1)}_{i;\mathrm{ad},{\mathrm{gls}}}(
\hat {A}_{i,{\mathrm{gls}}})&=&0,
\\
l_{\mathrm{RE}}^{(1)}(\hat{A}_{i,{\mathrm{ols}}})+\tilde{l}^{(1)}_{i;\mathrm{ad},{\mathrm{ols}}}(
\hat {A}_{i,{\mathrm{ols}}})&=&0,
\\
l_{\mathrm{RE}}^{(2)}(\hat{A})+\tilde{l}_{i;\mathrm{ad}}^{(2)}(
\hat{A})&<&0,
\end{eqnarray*}
where $\hat{A}\in\{\tilde{A}_{\mathrm{RE}}, \hat{A}_{i,{\mathrm{gls}}}, \hat{A}_{i,{\mathrm{ols}}}\}$
and $\tilde{A}_{\mathrm{RE}}$ is a solution to the REML estimation equation.
Hence, $\hat{A}_{\mathrm{RE}}$ is always larger than $\hat{A}_{i,{\mathrm{gls}}}$ or $\hat
{A}_{i,{\mathrm{gls}}}$ using the facts that $\hat{A}_{\mathrm{RE}}=\max\{0,\tilde
{A}_{\mathrm{RE}}\}$ and $\hat{A}_{i,{\mathrm{gls}}}$ or $\hat{A}_{i,{\mathrm{gls}}}$ are strictly
positive if $m>(4+p)/(1-q_i)$.

Finally, using that $0<\tilde{l}^{(1)}_{i;\mathrm{ad},{\mathrm{gls}}}\leq\tilde
{l}^{(1)}_{i;\mathrm{ad},{\mathrm{ols}}}$ for $A\geq0$, we have the result.
\end{pf*}

\begin{pf*}{Proof of Corollary to Theorem~\ref{th4}}
(i) Since for these two adjustment terms, $h_i(A)L_{\mathrm{RE}}(A) {|}_{A=0}=0$
and $h_i(A)L_{\mathrm{RE}}(A)\ge0$ for $A>0$, it suffices to show that $\lim_{A\rightarrow\infty} h_i(A)L_{\mathrm{RE}}(A)=0$.
For $h_i(A)$ given by (\ref{h.unbalance.ols}),
\[
(\ref{h.unbalance.ols})\leq (A+D_i)^2\Bigl(A+\sup
_{i\geq1}D_i\Bigr)^{({1}/{2})m q_i}\leq\Bigl(A+\sup
_{i\geq1}D_i\Bigr)^{2+({1}/{2})mq_i}.
\]

For (\ref{h.unbalance.gls}), we have
\begin{eqnarray*}
(\ref{h.unbalance.gls})&\leq& (A+D_i)^2\exp\biggl\{
\frac{1}{2}\int\Bigl(A+\inf_{i\geq1}D_i\Bigr)
q_i \operatorname{tr}\bigl(V^{-2}\bigr)\,dA\biggr\}
\\
&\leq&(A+D_i)^2\Bigl(A+\sup_{i\geq1}D_i
\Bigr)^{({1}/{2})m q_i} \exp\biggl[-\frac{m}{2}q_i\biggr]
\\
&&{}\times\exp\biggl[-\frac{1}{2}\inf_{i\geq1}D_i
q_i \operatorname{tr}\bigl(V^{-1}\bigr)\biggr]
\\
&\leq&\Bigl(A+\sup_{i\geq1}D_i\Bigr)^{2+({1}/{2})m q_i}.
\end{eqnarray*}
Using the fact $L_{\mathrm{RE}}(A)<C(A+\sup_{i\geq1}D_i)^{{p}/{2}}|
X^{\prime}X |^{-{1}/{2}}(A+\inf_{i\geq1}D_i)^{-{m}/{2}}$, we have
\[
0\leq h_i(A)L_{\mathrm{RE}}(A)\leq \Bigl(A+\sup
_{i\geq1}D_i\Bigr)^{2+({1}/{2})[m
q_i+p]}\Bigl(A+\inf
_{i\geq1}D_i\Bigr)^{-{m}/{2}}\bigl|X^{\prime}X\bigr|^{-{1}/{2}},
\]
so that, under mild regularity conditions,
\[
0\leq\lim_{A\rightarrow\infty} h_i(A)L_{\mathrm{RE}}(A)=\lim
_{A\rightarrow
\infty}A^{2+({1}/{2})[m q_i+p-m]}.
\]
Thus, if $2+\frac{1}{2}[m q_i+p-m]<0$, we have
\[
\lim_{A\rightarrow\infty} h_i(A)L_{\mathrm{RE}}(A)=0.
\]

We first show that $\hat{A}_{i;{\mathrm{gls}}}$ and $\hat{A}_{i;{\mathrm{gls}}}$ satisfy the
regularity conditions of Theorem~\ref{th3}.
Since $0<A<\infty$, we claim that $\tilde{l}_{i,\mathrm{ad}}^{k}(A)=O(1)\ (k=1,2,3)$, for large $m$, for both the GLS and OLS estimators of $\beta
$ using the following facts.

For the GLS estimator,
\begin{eqnarray*}
\tilde{l}_{i,\mathrm{ad}}^{(1)}(A)&=& \biggl(2-\frac{(1+z^2)}{4} \biggr)
\frac
{1}{A+D_i}+\frac{(1+z^2)}{4A}+\frac{1}{2}\operatorname{tr}\bigl[V^{-2}
\bigr]x_i^{\prime
}\bigl(X^{\prime}V^{-1}X
\bigr)^{-1}x_i,
\\
\tilde{l}_{i,\mathrm{ad}}^{(2)}(A)&=&- \biggl(2-\frac{(1+z^2)}{4}
\biggr)\frac
{1}{(A+D_i)^2}-\frac{(1+z^2)}{4A^2}
\\
&&{}-\operatorname{tr}\bigl[V^{-3}\bigr]x_i^{\prime}
\bigl(X^{\prime}V^{-1}X\bigr)^{-1}x_i\\
&&{}+
\frac
{1}{2}\operatorname{tr}\bigl[V^{-2}\bigr]x_i^{\prime}
\bigl(X^{\prime}V^{-1}X\bigr)^{-1}X^{\prime
}V^{-2}X
\bigl(X^{\prime}V^{-1}X\bigr)^{-1}x_i,
\\
\tilde{l}_{i,\mathrm{ad}}^{(3)}(A)&=& \biggl(2-\frac{(1+z^2)}{4} \biggr)
\frac
{2}{(A+D_i)^3}+\frac{(1+z^2)}{2A^3}
\\
&&{}+3\operatorname{tr}\bigl[V^{-4}\bigr]x_i^{\prime}
\bigl(X^{\prime
}V^{-1}X\bigr)^{-1}x_i\\
&&{}-2\operatorname{tr}
\bigl[V^{-3}\bigr]x_i^{\prime}\bigl(X^{\prime
}V^{-1}X
\bigr)^{-1}X^{\prime}V^{-2}X\bigl(X^{\prime}V^{-1}X
\bigr)^{-1}x_i
\\
&&{}\times\operatorname{tr}\bigl[V^{-2}\bigr] \bigl[x_i^{\prime}
\bigl(X^{\prime}V^{-1}X\bigr)^{-1}X^{\prime
}V^{-2}X
\bigl(X^{\prime}V^{-1}X\bigr)^{-1}\\
&&\hspace*{50pt}{}\times X^{\prime}V^{-2}X
\bigl(X^{\prime
}V^{-1}X\bigr)^{-1}x_i
\\
&&\hspace*{50pt}{}-x_i^{\prime}\bigl(X^{\prime}V^{-1}X
\bigr)^{-1}X^{\prime
}V^{-3}X\bigl(X^{\prime}V^{-1}X
\bigr)^{-1}x_i\bigr].
\end{eqnarray*}
For the OLS estimator,
\begin{eqnarray*}
\tilde{l}_{i,\mathrm{ad}}^{(1)}(A)&=& \biggl(2-\frac{(1+z^2)}{4} \biggr)
\frac
{1}{A+D_i}+\frac{(1+z^2)}{4A}
\\
&&{}+\frac{1}{2}\operatorname{tr}\bigl[V^{-2}\bigr]x_i^{\prime}
\bigl(X^{\prime}X\bigr)^{-1}X^{\prime
}VX
\bigl(X^{\prime}X\bigr)^{-1}x_i,
\\
\tilde{l}_{i,\mathrm{ad}}^{(2)}(A)&=&- \biggl(2-\frac{(1+z^2)}{4}
\biggr)\frac
{1}{(A+D_i)^2}-\frac{(1+z^2)}{4A^2}
\\
&&{}-\operatorname{tr}\bigl[V^{-3}\bigr]x_i^{\prime}
\bigl(X^{\prime}X\bigr)^{-1}X^{\prime}VX
\bigl(X^{\prime
}X\bigr)^{-1}x_i+\frac{1}{2}\operatorname{tr}
\bigl[V^{-2}\bigr] q_i,
\\
\tilde{l}_{i,\mathrm{ad}}^{(3)}(A)&=& \biggl(2-\frac{(1+z^2)}{4} \biggr)
\frac
{2}{(A+D_i)^3}+\frac{(1+z^2)}{2A^3}
\\
&&{}+3\operatorname{tr}\bigl[V^{-4}\bigr]x_i^{\prime}
\bigl(X^{\prime}X\bigr)^{-1}X^{\prime}VX
\bigl(X^{\prime
}X\bigr)^{-1}x_i-2\operatorname{tr}
\bigl[V^{-3}\bigr]q_i.
\end{eqnarray*}

In addition,
For GLS,
\[
\tilde{l}_{i,\mathrm{ad}}^{(4)}(A)=- \biggl(12-
\frac{3(1+z^2)}{2} \biggr)\frac
{1}{(A+D_i)^4}-\frac{3(1+z^2)}{2A^4}+
\tilde{l}_{3,i,\mathrm{ad},{\mathrm{gls}}}^{(4)}(A).
\]
For OLS,
\[
\tilde{l}_{i,\mathrm{ad}}^{(4)}(A)=-
\biggl(12-\frac{3(1+z^2)}{2} \biggr)\frac
{1}{(A+D_i)^4}-\frac{3(1+z^2)}{2A^4}+
\tilde{l}_{3,i,\mathrm{ad},\mathrm{ols}}^{(4)}(A),
\]
where
\begin{eqnarray*}
\tilde{l}_{3,i,\mathrm{ad},{\mathrm{gls}}}^{(4)}(A)&=&-12\operatorname{tr}\bigl[V^{-5}
\bigr]x_i^{\prime}\bigl(X^{\prime
}V^{-1}X
\bigr)^{-1}x_i
\\
&&{}+6\operatorname{tr}\bigl[V^{-3}\bigr] \bigl[ x_i^{\prime}
\bigl(X^{\prime}V^{-1}X\bigr)^{-1} X^{\prime}V^{-3}X
\bigl(X^{\prime}V^{-1}X\bigr)^{-1}x_i
\\
&&\hspace*{55pt}{}-x_i^{\prime}\bigl(X^{\prime}V^{-1}X
\bigr)^{-1} X^{\prime}V^{-2}X \bigl(X^{\prime
}V^{-1}X
\bigr)^{-1} X^{\prime}V^{-2}\\
&&\hspace*{178pt}{}\times X\bigl(X^{\prime}V^{-1}X
\bigr)^{-1}x_i\bigr]
\\
&&{}+9\operatorname{tr}\bigl[V^{-4}\bigr] x_i^{\prime}
\bigl(X^{\prime}V^{-1}X\bigr)^{-1} X^{\prime}V^{-2}X
\bigl(X^{\prime}V^{-1}X\bigr)^{-1}x_i
\\
&&{}+\operatorname{tr}\bigl[V^{-2}\bigr] \bigl[3x_i^{\prime}
\bigl(X^{\prime}V^{-1}X\bigr)^{-1} X^{\prime}V^{-4}X
\bigl(X^{\prime}V^{-1}X\bigr)^{-1}x_i
\\
&&\hspace*{48pt}{}-4x_i^{\prime}\bigl(X^{\prime}V^{-1}X
\bigr)^{-1} X^{\prime}V^{-2}X \bigl(X^{\prime
}V^{-1}X
\bigr)^{-1} X^{\prime}V^{-3}\\
&&\hspace*{48pt}{}\times X \bigl(X^{\prime}V^{-1}X
\bigr)^{-1}x_i
\\
&&\hspace*{48pt}{}-4x_i^{\prime}\bigl(X^{\prime}V^{-1}X
\bigr)^{-1} X^{\prime}V^{-3}X \bigl(X^{\prime
}V^{-1}X
\bigr)^{-1} X^{\prime}V^{-2}\\
&&\hspace*{48pt}{}\times X \bigl(X^{\prime}V^{-1}X
\bigr)^{-1}x_i
\\
&&\hspace*{48pt}{}+3x_i^{\prime}\bigl(X^{\prime}V^{-1}X
\bigr)^{-1} X^{\prime}V^{-2}X \bigl(X^{\prime
}V^{-1}X
\bigr)^{-1} X^{\prime}V^{-2}\\
&&\hspace*{48pt}{}\times X \bigl(X^{\prime}V^{-1}X
\bigr)^{-1}
\\
&&\hspace*{146pt}{}\times X^{\prime}V^{-2}X \bigl(X^{\prime}V^{-1}X
\bigr)^{-1} x_i\bigr],
\\
\tilde{l}_{3,i,\mathrm{ad},{\mathrm{ols}}}^{(4)}(A)&=&-12\operatorname{tr}\bigl[V^{-5}
\bigr]x_i^{\prime}\bigl(X^{\prime
}X\bigr)^{-1}X^{\prime}VX
\bigl(X^{\prime}X\bigr)^{-1}x_i
\\
&&{}+9\operatorname{tr}\bigl[V^{-4}\bigr]q_i.
\end{eqnarray*}

Using the above facts, we can prove that $|\tilde{l}_{i,\mathrm{ad},{\mathrm{gls}}}^{(j)}|$
and $|\tilde{l}_{i,\mathrm{ad},{\mathrm{ols}}}^{(j)}|$ are bounded for $j=1,2,3,4$ under
the regularity conditions R2--R4. Similarly, we can show that the $g$th
powers of
$\sup_{A/2<{A}^*<2A}\frac{1}{m}|\tilde{l}_{i,\mathrm{ad};h}^{(5)}({A}^*)|$ with
$h={\mathrm{gls}}, {\mathrm{ols}}$ are bounded for any fixed $g>0$. Thus, the new area specific
adjustment terms satisfy the regularity condition R1. Thus, an
application of Theorem~\ref{th4} leads to (ii)--(iv) of the
Corollary to Theorem~\ref{th4}.
\end{pf*}

\section{Proof of Lemma~\texorpdfstring{\lowercase{\protect\ref{le1}}}{1}}

The proof of (\ref{u1.2}) is much more complex due to the
dependence of $\hat{A}$ and $y_i$. We use the following lemma
repeatedly for proving (\ref{u1.2}). For a proof of Lemma~\ref{ST1},
see \citet{ST}.
%
\begin{lem}
\label{ST1}
Let $Z\sim N(0,\Sigma)$. Then for symmetric matrices $Q$, $U$ and $W$,
\begin{eqnarray*}
E\bigl[\bigl(Z^{\prime}QZ\bigr) \bigl(Z^{\prime}UZ\bigr)
\bigr]&=&2\operatorname{tr}(Q\Sigma U\Sigma) +\operatorname{tr}(Q\Sigma )\operatorname{tr}(U\Sigma),
\\
E\bigl[\bigl(Z^{\prime}QZ\bigr) \bigl(Z^{\prime}UZ\bigr)
\bigl(Z^{\prime}WZ\bigr)\bigr] &=&8\operatorname{tr}(Q\Sigma U\Sigma W\Sigma)\\
&&{}+2\bigl\{\operatorname{tr}(Q
\Sigma U\Sigma) \operatorname{tr}(W\Sigma)+ \operatorname{tr}(Q\Sigma W\Sigma) \operatorname{tr}(U\Sigma)
\\
&&\hspace*{110pt}{}+\operatorname{tr}(U\Sigma W\Sigma) \operatorname{tr}(Q\Sigma)\bigr\} \\
&&{}+\operatorname{tr}(Q\Sigma) \operatorname{tr}(U\Sigma) \operatorname{tr}(W\Sigma).
\end{eqnarray*}
\end{lem}

The proof also needs the following lemma, which is immediate from
Theorem~2.1 of Das, Jiang and Rao (\citeyear{Das}).
%
\begin{lem}
\label{Das1}
Assume the following regularity conditions:
\begin{longlist}[1.]
\item[1.]$\tilde{l}_{i,\mathrm{ad}}(A)$, which is free of $y$, is four times
continuously differentiable with respect to $A$,
\item[2.] the $g$th power of the following are bounded: $\frac{1}{\sqrt
{m}}|\tilde{l}_{i,\mathrm{ad}}^{(1)}(A)|$, $\frac{1}{m}|\tilde
{l}_{i,\mathrm{ad}}^{(2)}(A)|$, $\frac{1}{m}|\tilde{l}_{i,\mathrm{ad}}^{(3)}(A)|$,
and $\frac{1}{m}\sup_{A/2<\tilde{A}<2A}\llvert \tilde
{l}^{(4)}_{i,\mathrm{ad}}(A) |_{A=\tilde{A}}\rrvert $ (fixed $g>0$),
\item[3.]$A\in\Theta_0$, the interior of $\Theta$, that is, $0<A<\infty$.
\end{longlist}
Then:

\textup{(i)} there is $\hat{A}_i$ such that for any $0<\rho<1$, there is a set
$\Lambda$ satisfying for large $m$ and on $\Lambda$,
$\hat{A}\in\Theta, l^{(1)}(A) |_{\hat{A}}=0$, $\sqrt{m}|\hat
{A}_i-A|<m^{{(1-\rho)}/{2}}$, and
\[
\hat{A}_i-A=I+\mathit{II}+\mathit{III}+r,
\]
where $I=-E[l^{(2)}]^{-1}l^{(1)}$,
$\mathit{II}=E[l^{(2)}]^{-2}l^{(2)}l^{(1)}-E[l^{(2)}]^{-1}l^{(1)}$,
$\mathit{III}=\break  -\frac{1}{2}E[l^{(2)}]^{-3}\{l^{(1)}\}^2l^{(3)}$, and $r\leq
m^{-3\rho/2}u$ with $E[|u|^g]$ bounded;

\textup{(ii)} $P(\Lambda^{c})\leq m^{-\tau/2 g}C$, where $\tau=1/4 \wedge(1-\rho
)$.
\end{lem}

First note that
\[
E\bigl[\hat{u}_{1i}^2\bigr]=m\sigma_i^{-2}
\biggl\{B_i^2 T_1+2B_i
\frac
{D_i}{(A+D_i)^2}T_2+\frac{D_i^2}{(A+D_i)^4}T_3 \biggr\},
\]
where $T_1=E[x_{i}^{\prime}(\hat{\beta}-\beta)^2]$, $T_2= E[(\hat
{A}-A)x_{i}^{\prime}(\hat{\beta}-\beta)(y_i-x_{i}^{\prime}\beta)]$ and
$T_3= E[(\hat{A}-A)^2(y_i-x_i^{\prime} \beta)^2]$. We now simplify
these three terms.

We first prove that
%
\begin{equation}
E[T_1]=x^{\prime}_i\operatorname{Var}(\tilde{
\beta})x_{i}+O\bigl(m^{-2}\bigr),\label{T.1}
\end{equation}
where $\operatorname{Var}(\tilde{\beta})=(X^{\prime}X)^{-1}X^{\prime}VX(X^{\prime
}X)^{-1}$ if $\tilde{\beta}$ is the OLS estimator of $\beta$ and
$(X^{\prime}V^{-1}X)^{-1}$ if $\tilde{\beta}$ is the GLS estimator of
$\beta$.

Note that
\begin{eqnarray*}
E\bigl[\bigl\{x_i^{\prime}(\hat{\beta}-{\beta})\bigr
\}^2\bigr]&=&E\bigl[\bigl\{x_i^{\prime}(\tilde {
\beta}-{\beta})\bigr\}^2\bigr]+E\bigl[\bigl\{x_i^{\prime}(
\hat{\beta}-\tilde{\beta})\bigr\} ^2\bigr]
\\
&=&x^{\prime}_i\operatorname{Var}(\tilde{\beta})x_{i}+E\bigl[\bigl
\{x_i^{\prime}(\hat{\beta }-\tilde{\beta})\bigr\}^2
\bigr],
\end{eqnarray*}
and we have the following facts:
%
\begin{equation}
E\bigl[\bigl\{x_i^{\prime}\bigl(\hat{\beta}(
\hat{A}_{1},\ldots,\hat{A}_{m})-\tilde {\beta}(A)\bigr)
\bigr\}^2\bigr]\leq E\bigl[\bigl\{x_i^{\prime}
\bigl(\hat{\beta}(\hat{A}_{U})-\tilde {\beta}\bigr)\bigr
\}^2\bigr],\label{xbeta.2}
\end{equation}
where $\hat{A}_{U}=\argmax_{\hat{A}_i} |x_i^{\prime}(\hat{\beta}(\hat
{A}_{1},\ldots,\hat{A}_{m})-\tilde{\beta}(A))|$.

We have $\frac{\partial\tilde{\beta}}{\partial A}=H(y-X\beta)$, where
$H=0$ for the OLS estimator of $\beta$ and
\[
H=\bigl(X^{\prime}V^{-1}X\bigr)^{-1}
X^{\prime}V^{-2}X\bigl( X^{\prime
}V^{-1}X
\bigr)^{-1}X^{\prime}V^{-1}-\bigl(X^{\prime}V^{-1}X
\bigr)^{-1} X^{\prime}V^{-2},
\]
the GLS estimators of $\beta$.

Using the Taylor series expansion, we have
%
\begin{equation}
x_{i}^{\prime}\bigl(\hat{\beta}(\hat{A}_{U})-
\tilde{\beta}\bigr)=(\hat {A}_{U}-A)x_{i}^{\prime}Hy+r_1,\label{x.beta.app}
\end{equation}
where
$|r_1|=\frac{1}{2}(\hat{A}_U-A)^2 x_i^{\prime}\frac{\partial
H}{\partial A} |_{A=A^{*}}y$ with $A^{*}\in(A,\hat{A}_U)$ and
\begin{eqnarray*}
\frac{\partial H}{\partial A}&=&2\bigl(X^{\prime}V^{-1}X\bigr)^{-1}X^{\prime
}V^{-2}
\bigl(X\bigl(X^{\prime}V^{-1}X\bigr)^{-1}X^{\prime} V^{-1}-I
\bigr)\\
&&{}\times V^{-1}\bigl(X\bigl(X^{\prime}V^{-1}X
\bigr)^{-1}X^{\prime}V^{-1}-I\bigr).
\end{eqnarray*}

Let $H_s^{(1)}$ be the matrix with $(i,j)$ components given by
\[
\sup_{A/2<A^*<2A} \biggl\{\frac{\partial H}{\partial A} \bigg|_{A=A^*}
\biggr\}_{(i,j)},
\]
where ${Q}_{(i,j)}$ is $(i,j)$ component of a matrix $Q$.
Under the regularity conditions~R3--R4, we can show that the
components of $H_s^{(1)}$ are bounded and of order $O(m^{-1})$ using an
argument similar to that given in Proposition~3.2 of Das, Jiang and Rao (\citeyear{Das}).
Using the facts that $HX=0$, $x_{i}^{\prime}HVH^{\prime}x_i=O(m^{-1})$,
we have
\begin{eqnarray*}
E\bigl[\bigl\{x_{i}^{\prime}(\hat{\beta}-\tilde{\beta})\bigr
\}^2\bigr] &\leq&E \bigl[(\hat{A}_{U}-A)^2
\bigl(x^{\prime}_i Hy\bigr) \bigl(y^{\prime}H^{\prime
}x_i
\bigr) \bigr]\\
&&{}+ 2E \bigl[(\hat{A}_U-A)^3
\bigl(x^{\prime}_i Hy\bigr) \bigl(y^{\prime}
\bigl[H_s^{(1)}\bigr]^{\prime
}x_i \bigr)
\bigr]
\\
&&{}+E \biggl[(\hat{A}_U-A)^4 \biggl(x^{\prime}_i
\frac{\partial H}{\partial
A}y \biggr) \bigl(y^{\prime}\bigl[H_s^{(1)}
\bigr]^{\prime}x_i \bigr) \biggr]
\\
&\leq&E\bigl[(\hat{A}_U-A)^2\bigr]x_{i}^{\prime}HVH^{\prime}x_i+E
\bigl[|\hat {A}_U-A|^3\bigr]x_{i}^{\prime}HV
\bigl[H_s^{(1)}\bigr]^{\prime}x_i
\\
&&{}+E\bigl[(\hat{A}_U-A)^4\bigr]x_{i}^{\prime}H_s^{(1)}V
\bigl[H_s^{(1)}\bigr]^{\prime}x_i,
\\
&=&O\bigl(m^{-2}\bigr).
\end{eqnarray*}

Thus, this completes the proof of (\ref{T.1}).

Next, we simplify $E[T_2]$. Let $l_{i;\mathrm{ad}}$ denote the adjusted residual
log-likelihood function. Then $l_{i;\mathrm{ad}}=l_{\mathrm{RE}}+\tilde{l}_{i;\mathrm{ad}}$, where
$l_{\mathrm{RE}}$ is the residual log-likelihood function and $\tilde
{l}_{i;\mathrm{ad}}=\log h_i(A)$. Define $I_{F}=-1/E[\frac{\partial^2
l}{\partial A^2}]$. For notational simplicity, we set $l_{i;\mathrm{ad}}\equiv
l_{\mathrm{ad}}$ and $\tilde{l}_{i;\mathrm{ad}}\equiv\tilde{l}_{\mathrm{ad}}$. Since $\tilde
{l}_{\mathrm{ad}}$ is bounded and free from $y$, we obtain the following using
Lemma~\ref{Das1},
\[
\hat{A}-A=\frac{\partial l_{\mathrm{ad}}}{\partial
A}I_{F}+r_{2.1}=l^{(1)}_{\mathrm{RE}}I_{F}+r_{2.2},
\]
where $l^{(1)}_{\mathrm{RE}}=\frac{\partial l_{\mathrm{RE}}}{\partial A}=\frac
{1}{2}[y^{\prime}P^2y-\operatorname{tr}(P)]$ and $E[|r_{2.2}|]=O(m^{-1})$ when $\rho$
is taken as $3/4$ in Lemma~\ref{Das1}.

Since $\hat{A}$ is translation invariant and even function, we can
substitute $\hat{A}(Z)-A$ for $\hat{A}(y)-A$, where $Z=y-X \beta\sim N(0,V)$.
Thus,
\begin{eqnarray*}
x_i^{\prime}(\hat{\beta}-\beta)&=&x^{\prime}_i
\bigl(X^{\prime}\hat {V}^{-1}X\bigr)^{-1}X^{\prime}
\hat{V}^{-1}Z
\\
&=&\lambda_i^{\prime}X\bigl(X^{\prime}
\hat{V}^{-1}X\bigr)^{-1}X^{\prime}\hat
{V}^{-1}Z
\\
&=&\lambda_i^{\prime}X\bigl(X^{\prime}{V}^{-1}X
\bigr)^{-1}X^{\prime} {V}^{-1}Z+r_{1.2}Z,
\end{eqnarray*}
where $\lambda_i$ denotes a $m\times1$ vector with $i$ component 1 and
the rest 0 and
$r_{1.2}Z\leq(\hat{A}_U-A)x_{i}^{\prime}HZ+r_1$.

Hence,
\begin{eqnarray*}
E[T_2]&\leq &E \bigl[\bigl(l^{(1)}_{\mathrm{RE}}
I_{F}+r_{2.2}\bigr) \bigl\{\lambda _i^{\prime}
X \bigl(X^{\prime}V^{-1}X\bigr)^{-1}X^{\prime}V^{-1}
Z+r_{1.2}Z \bigr\}\bigl(\lambda_i^{\prime} Z\bigr)
\bigr]
\\
&=&I_F \bigl\{ E \bigl[l^{(1)}_{\mathrm{RE}}Z^{\prime}E_i
X\bigl(X^{\prime
}V^{-1}X\bigr)^{-1}X^{\prime}V^{-1}Z
\bigr] \bigr\}+ E \bigl[(\hat {A}_U-A)r_{1.2}Z\bigl(
\lambda_i^{\prime} Z\bigr) \bigr]
\\
&&{}+E \bigl[r_{2.2}Z^{\prime}E_i X
\bigl(X^{\prime}V^{-1}X\bigr)^{-1}X^{\prime
}V^{-1}
Z \bigr]
\\
&=&I_F T_{2.1}+T_{2.2}+T_{2.3},
\end{eqnarray*}
where $E_i$ denotes a $m\times m$ matrix with the $(i,i)$ component one
and rest zeroes.

Using Lemma~\ref{ST1} and the following facts:
\begin{longlist}[(ii)]
\item[(i)] $PVP=P$,
\item[(ii)] $\operatorname{tr}[C_iV]$ and $\operatorname{tr}[P^{2}VC_iV]$ are of order $O(m^{-1})$,
under the regularity conditions,
\end{longlist}
we have
\begin{eqnarray*}
T_{2.1}&=&\tfrac{1}{2} \bigl\{E\bigl[\bigl(Z^{\prime}P^2Z
\bigr) \bigl(Z^{\prime}C_i Z\bigr)\bigr]-\operatorname{tr}[P]E
\bigl[Z^{\prime}C_i Z\bigr] \bigr\}
\\
&=&\operatorname{tr}\bigl[P^2VC_iV\bigr]+\tfrac{1}{2}\operatorname{tr}
\bigl[P^2V\bigr]\operatorname{tr}[C_iV]-\tfrac{1}{2}\operatorname{tr}[P]\operatorname{tr}[C_iV]
\\
&=&O\bigl(m^{-1}\bigr),
\end{eqnarray*}
where $C_i=E_i X(X^{\prime}V^{-1}X)^{-1}X^{\prime}V^{-1}$.

Using $\lambda^{\prime}_i\frac{\partial H}{\partial A}=O(m^{-1})$, we have
\begin{eqnarray*}
T_{2.2}&=&E\bigl[(\hat{A}_U-A)r_{1.2}Z \bigl(
\lambda_i^{\prime}Z\bigr)\bigr]
\\
&=&E \bigl[(\hat{A}_U-A)^2 \bigl(\lambda^{\prime}_i
XHZ\bigr) \bigl( \lambda_i^{\prime
}Z\bigr) \bigr]+E \bigl[(
\hat{A}_U-A)r_{1}\bigl( \lambda_i^{\prime}Z
\bigr) \bigr]
\\
&\leq& E\bigl[(\hat{A}_U-A)^2\bigr]E
\bigl[Z^{\prime}H^{\prime}X^{\prime}E_iZ\bigr]+ E
\biggl[(\hat{A}_U-A)^3 \biggl(\lambda_i^{\prime}X
\frac{\partial
H}{\partial A} \bigg|_{A=A}Z \biggr) \bigl(\lambda_{i}^{\prime}Z
\bigr) \biggr]
\\
&=&O\bigl(m^{-2}\bigr).
\end{eqnarray*}

Using $E[|r_{2.2}|]=O(m^{-1})$,
\[
T_{2.3}=E\bigl[r_{2.2}Z^{\prime}C_i Z
\bigr]\leq E\bigl[|r_{2.2}|\bigr]\operatorname{tr}[C_i V]=O\bigl(m^{-2}
\bigr).
\]
Therefore,
\[
E[T_{2}]\leq O\bigl(m^{-2}\bigr).
\]

Hence, using the above results and $E[T_{2}]\geq O(m^{-2})$ with same
calculation, we have
%
\begin{equation}
E[T_2]=O\bigl(m^{-2}\bigr).\label{T.2}
\end{equation}

Since $I_F$ is of order $O(m^{-1})$, we have
\begin{eqnarray*}
E[T_3]&=&E\bigl[(\hat{A}-A)^2\bigl(y_i-x_i^{\prime}
\beta\bigr)^2\bigr]
\\
&=&E\bigl[\bigl(I_F l_{\mathrm{RE}}^{(1)}+r_{2.2}
\bigr)^2\lambda_i^{\prime}ZZ^{\prime}
\lambda_i \bigr]
\\
&=&I_{F}^2 \bigl\{ \tfrac{1}{4}E\bigl[
\bigl(Z^{\prime}P^2Z\bigr) \bigl(Z^{\prime}P^2Z
\bigr) \bigl(Z^{\prime
}E_iZ\bigr)\bigr]-\tfrac{1}{2}E
\bigl[\bigl(Z^{\prime}P^2Z\bigr) \bigl(Z^{\prime}E_iZ
\bigr)\bigr]\operatorname{tr}[P]
\\
&&\hspace*{205pt}{}+\tfrac{1}{4}E\bigl[Z^{\prime}E_iZ
\bigr]\operatorname{tr}[P]^2 \bigr\}\\
&&{}+I_F E\bigl[r_{2.2}
\bigl(Z^{\prime
}P^2Z-\operatorname{tr}[P]\bigr)Z^{\prime}E_iZ
\bigr]+E\bigl[r_{2.2}^2 Z^{\prime}E_iZ
\bigr]
\\
&\leq&I_F^2 \Gamma+I_F E[|r_{2.2}|]
\bigl\{ 2\operatorname{tr}\bigl[P^2VE_iV\bigr]+\operatorname{tr}\bigl[P^2V
\bigr]\operatorname{tr}[E_iV]-\operatorname{tr}[P]\operatorname{tr}[E_iV]\bigr\}
\\
&&{}+E
\bigl[r_{2.2}^2\bigr]\operatorname{tr}[E_iV]
\\
&=&\Gamma I_{F}^2+O\bigl(m^{-2}\bigr).
\end{eqnarray*}

Using Lemma~\ref{ST1} and the following facts:

\begin{longlist}[(iii)]
\item[(i)] $PVP=P$,
\item[(ii)] $\operatorname{tr}(E_iV)=(A+D_i)$, and
\item[(iii)] $|\operatorname{tr}(P^k)-\operatorname{tr}(V^{-k})|=O(1)$, for $k\geq1$,\vadjust{\goodbreak}
\end{longlist}
we have
\begin{eqnarray*}
\Gamma&=& \bigl\{ \tfrac{1}{4}E\bigl[\bigl(Z^{\prime}P^2Z
\bigr) \bigl(Z^{\prime}P^2Z\bigr) \bigl(Z^{\prime
}E_iZ
\bigr)\bigr]-\tfrac{1}{2}E\bigl[\bigl(Z^{\prime}P^2Z\bigr)
\bigl(Z^{\prime}E_iZ\bigr)\bigr]\operatorname{tr}[P] \\
&&\hspace*{189pt}{}+\tfrac{1}{4}E
\bigl[Z^{\prime}E_iZ\bigr]\operatorname{tr}[P]^2 \bigr\}
\\
&=&\tfrac{1}{4} \bigl[8\operatorname{tr}\bigl(P^2VP^2E_iV
\bigr)+2\bigl\{ \operatorname{tr}\bigl(P^2VP^2V\bigr)\operatorname{tr}(E_iV)+2\operatorname{tr}
\bigl(P^2VE_iV\bigr)\operatorname{tr}\bigl(P^2V\bigr)\bigr
\}\\
&&\hspace*{235pt}{}+\operatorname{tr}\bigl(P^2V\bigr)^2\operatorname{tr}(E_iV) \bigr]
\\
&&{}-\operatorname{tr}\bigl(P^2VE_iV\bigr)\operatorname{tr}(P) -\tfrac{1}{2}\operatorname{tr}
\bigl(P^2V\bigr)\operatorname{tr}(E_iV)\operatorname{tr}(P)+\tfrac
{1}{4}\operatorname{tr}(P)^2\operatorname{tr}(E_iV)
\\
&=&2\operatorname{tr}\bigl(P^3VE_iV\bigr)+\tfrac{1}{2}\operatorname{tr}
\bigl(P^2\bigr)\operatorname{tr}(E_iV)=\tfrac
{1}{2}\operatorname{tr}
\bigl(P^2\bigr)\operatorname{tr}(E_iV)+O(1)\\
&=&\tfrac{1}{2}\operatorname{tr}
\bigl(V^{-2}\bigr) (A+D_i)+O(1).
\end{eqnarray*}

Hence,
%
\begin{equation}\qquad
E[T_3]=I_{F}^2\frac{1}{2}\operatorname{tr}
\bigl(V^{-2}\bigr) (A+D_i)+O\bigl(m^{-2}\bigr)=
\frac
{2(A+D_i)}{\operatorname{tr}(V^{-2})}+O\bigl(m^{-2}\bigr).\label{T.3}
\end{equation}

Thus, we can show (\ref{u1.2}) using (\ref{T.1}), (\ref{T.2}) and (\ref
{T.3}).
\end{appendix}

\section*{Acknowledgments}
M. Yoshimori conducted this research while visiting the
University of Maryland, College Park, USA, as a research scholar under
the supervison of the second author.
The authors thank Professor Yutaka Kano for reading an earlier draft of
the paper and making constructive comments. We are also grateful to
three referees, an Associate Editor and Professor Runze Li, Coeditor,
for making a number of constructive suggestions, which led to a
significant improvement of our paper.

\begin{supplement}[id=suppA]
\stitle{Supplemental proof}
\slink[doi]{10.1214/14-AOS1219SUPP} 
\sdatatype{.pdf}
\sfilename{aos1219\_supp.pdf}
\sdescription{We provide a proof of Theorem~\ref{th4}.}
\end{supplement}

%


\printaddresses
\end{document}